\documentclass[preprint]{elsarticle}

\usepackage{hyperref}
\usepackage{amsmath}
\usepackage{verbatim}
\usepackage{amssymb}

\usepackage[usenames, dvipsnames]{color}
\usepackage{mathtools}

\newcommand{\aside}[1]{\textcolor{black}{#1}}

\newcommand{\R}{\mathbb{R}}

\newcommand{\beq}{\begin{equation}}
\newcommand{\eeq}{\end{equation}}
\newcommand{\bseq}{\begin{equation*}}
\newcommand{\eseq}{\end{equation*}}
\newcommand{\bln}{}
\newcommand{\eln}{}
\newcommand{\bal}{\begin{aligned}}
\newcommand{\eal}{\end{aligned}}
\newcommand{\bgat}{\begin{gathered}}
\newcommand{\egat}{\end{gathered}}

\renewcommand{\bar}{\overline}

\journal{Journal of Computational Physics}

\begin{document}

\begin{frontmatter}

\title{A Kernel-based Lagrangian Method for Imperfectly-mixed Chemical Reactions \tnoteref{mytitlenote}}
\tnotetext[mytitlenote]{The authors were supported by the National Science Foundation under awards EAR-1417145, DMS-1211667, and DMS-1614586.}

\author{Michael J. Schmidt\fnref{hydro,ams}}
\ead{mschmidt1@mines.edu}
\author{Stephen Pankavich\fnref{ams}}
\ead{pankavic@mines.edu}
\author{David A. Benson\fnref{hydro}}
\ead{dbenson@mines.edu}
\address{Colorado School of Mines\\ 1500 Illinois St.\\ Golden, CO 80401}
\fntext[hydro]{Hydrologic Science and Engineering Program, Department of Geology and Geological Engineering, Colorado School of Mines, Golden, CO, 80401, USA}
\fntext[ams]{Department of Applied Mathematics and Statistics, Colorado School of Mines, Golden, CO, 80401, USA}

\begin{abstract}

Current Lagrangian (particle-tracking) algorithms used to simulate diffusion-reaction equations must employ a certain number of particles to properly emulate the system dynamics---particularly for imperfectly-mixed systems.  The number of particles is tied to the statistics of the initial concentration fields of the system at hand.  Systems with shorter-range correlation and/or smaller concentration variance require more particles, potentially limiting the computational feasibility of the method. For the well-known problem of bimolecular reaction, we show that using kernel-based, rather than Dirac delta, particles can significantly reduce the required number of particles.  We derive the fixed width of a Gaussian kernel for a given reduced number of particles that analytically eliminates the error between kernel and Dirac solutions at any specified time.  We also show how to solve for the fixed kernel size by minimizing the squared differences between solutions over any given time interval. Numerical results show that the width of the kernel should be kept below about 12\% of the domain size, and that the analytic equations used to derive kernel width suffer significantly from the neglect of higher-order moments.  The simulations with a kernel width given by least squares minimization perform better than those made to match at one specific time.  A heuristic time-variable kernel size, based on the previous results, performs on par with the least squares fixed kernel size.

\end{abstract}


\end{frontmatter}


\section{Introduction}
Chemical reactions can be simulated using several methods.  In particular, the (well-mixed) thermodynamic governing equations can be approximated using classical grid-based, Eulerian methods such as finite-differences or finite-elements. Below the scale of discretization---the subgrid scale---these techniques assume perfect mixing because at each node (or element, or cell), the concentration of each species is represented by a single, average, number. The effects of the imperfect mixing of concentration fluctuations cannot be resolved at the subgrid level except by empirical adjustments to the thermodynamic reaction rates \cite[e.g.,][]{battiato1,battiato2,Schwede_sample_pdf, Chiogna2013}. Grid-based methods must also deal with the difficulties of the hyperbolic portion of transport when reacting fluids are in motion \cite[e.g.,][]{Sweby1984,Leonard1991,Leveque2002,Bokanowski2007}. 

A number of particle-based transport and reaction algorithms have been developed to deal with imperfect mixing and high-P\'eclet number flows. These particle reaction algorithms are largely extensions of the Gillespie algorithm \cite{Gillespie1977,Gillespie2013}, with adjustments to handle spatial fluctuations of concentration in different ways.   These include particle-on-lattice methods that share some of the subgrid averaging effects as the Eulerian methods and have been shown to exhibit further inaccuracies when the lattice is refined below some minimum size \cite{Isaacson, HellanderPRE}. A more accurate algorithm is gridless and calculates the exact probabilities of particle collision and reaction.  The collisions are treated as independent \cite{Benson_react} or biased by nearest-neighbor status \cite{vanZon}.  Because collisions require a physically-based random motion, this method can be applied at any scale and under any type of diffusion \cite{Bolster2012}.  However, particle methods are not without their limitations.  If the reactions are simulated as a strict birth/death process, then the resolution for lower concentrations is limited to $\mathcal{O}(1/N)$, where $N$ is the particle number.  This  can be addressed by re-casting the birth/death process as a probabilistic, real-valued mass-reduction of individual particles \cite{Bolster_mass}.  Furthermore, reactions that are more complex than bi-molecular seemingly require a calculation of the collision of multiple particles, which can be very time-consuming. This problem can be alleviated by placing multiple chemical species on each particle and limiting the particle collision process to mass exchange and local mixing \cite{Benson_arbitrary}. 

One remaining problem with particle methods is that the initial particle number is not necessarily a modeling choice.  {\em Paster et al.} \cite{Paster_JCP} show that the required particle density in a finite domain is inversely proportional to the variance and the correlation length of a species' concentration fluctuations.  In short, ``smoother'' fields defined by lower variance and short-range correlation require more particles---if they are treated as Dirac delta functions.  Uniform (constant) initial concentrations would require an infinite number of particles, so clearly the method requires an alternative way to treat particles and their interaction.  An interesting extension treats each particle as a kernel that possesses some mean position, mass, and shape, imparting some degree of smoothness to each particle \cite{Dani_kernel}.  It remains to be shown, however, how the particle number changes with the shape of the kernel, or more precisely: how closely can two simulations be made to match when one has a given number of Dirac particles and another has a smaller number of kernel-smoothed particles?  

Thus, we intend to precisely address this questions within the current article. First, in Section \ref{sec:analytic} we define an initial value problem for diffusion-driven bimolecular reactions that admits analytic solutions and has been numerically approximated by many methods \cite{kang_prl,Benson_react,Bolster_mass}. We choose a problem that very clearly displays the reduced reaction rates engendered by imperfect mixing \cite{Ovchinnikov1978,KANG1984,Toussaint}.  In Section \ref{sec:analytic} we also define the ``match'' between two different solutions. In Section \ref{sec:particle} we define the numerical solutions for the Dirac and Gaussian kernel particles and identify the shape of the kernel based on user-defined particle numbers. In Section \ref{sec:properror} we show that our approach and methods are well-founded by proving a quantitative bound on the error in our approximation methods.  Finally, in Sections \ref{sec:results} and \ref{sec:concl} we show the results of ensembles of particle solutions and the approximate analytic equations for concentrations given Dirac versus Gaussian kernels of various size and provide some caveats for kernel use.
 
\section{Analytic model and methods}\label{sec:analytic}

\subsection{Governing equations}
We consider a bimolecular reactive system where the transport of constituents is driven by diffusion. Two components in
this system, $A$ and $B$, react kinetically and irreversibly with one another such that $A + B \to \emptyset$.
For an infinite, $d$-dimensional space, the governing equation for transport of the components is the diffusion-reaction equation (DRE)
\beq
\label{DRE}
\frac{\partial C_i}{\partial t} - D\Delta C_i = - k C_A C_B, \quad i = A, B,  \quad x\in\R^d, \quad t > 0,
\eeq
where the concentration of species $i$ is given by $C_i = C_i(t,x)$ [mol L$^{-d}$], $D$ is the diffusion coefficient [L$^2$ T$^{-1}$], and $k$ [L$^d$ mol$^{-1}$ T$^{-1}$] is the reaction rate constant. The reaction rate $r = kC_AC_B$ is assumed to be given by the law of mass action. \aside{If $C_0$ is the initial mean concentration of both $A$ and $B$ particles, this system can be characterized using a single dimensionless quantity, known as the Damk\"ohler number, $Da:=kC_0a^2/D$, where $a$ is some characteristic length-scale of diffusive mixing, typically given by the concentration autocorrelation length.}

\subsection{Exact solution for well-mixed system}
When the system is well-mixed in space at all times, concentrations are spatially homogeneous (i.e., constant in $x$), and the diffusion term in (\ref{DRE}) can be dropped. For this case, called the thermodynamic limit, well-known analytical solutions exist. For the simplest case, in which the initial average concentrations are equal $C_A(0,x) = C_B (0,x) = C_0 > 0$, Equation (\ref{DRE}) becomes $\frac{d C_i}{d t} = -k C_i^2$ , and the exact solution is merely
\beq
\label{wellmixed}
C_i(t) = \frac{C_0}{1+ C_0kt},
\eeq
which tends to zero like $t^{-1}$ for large time.

\subsection{Solution for a system with concentration fluctuations}
Keeping the thermodynamic limit in mind, we are interested in the more complex case where concentrations have some
random initial distribution, and the system is not continuously well-mixed. For this setup, concentrations evolve by diffusion
as described by the DRE (1). We may decompose concentrations into mean and fluctuation terms, such that
\beq
\label{fluctuation}
C_i(t,x) = \overline{C}_i(t) + C'_i(t,x),
\eeq
where the overbar refers to the ensemble average and the prime to the zero-mean fluctuations about the average. We
assume the system is ergodic, such that the ensemble average and spatial average are interchangeable, and only depend on $t$.
Restrictions associated with the assumption of ergodicity are discussed, for example, by \textit{Tartakovsky et al.} \cite{wrr_incomplete}. Broadly speaking, the larger the physical domain, the better the assumption of ergodicity, as larger domains can allow for a full ensemble of statistics. 
We shall also assume that the initial concentration distributions of species $A$ and $B$ have identical spatial autocovariance at $t = 0$. 
These assumptions allow one to obtain an analytical solution for the problem but are not required for numerical simulations. With these assumptions, the spatial averages of the initial concentrations of $A$ and $B$ are equal, i.e.,
\beq
\overline{C}_A(0,x) = \overline{C}_B (0,x) \equiv C_0,
\eeq
and clearly, due to the 1:1 stoichiometric ratio of the reaction, their average concentrations are equal for all $t \geq 0$, i.e.,
\beq
\overline{C}_A(t) = \overline{C}_B (t) = \overline{C}(t).
\eeq
Substituting (\ref{fluctuation}) into (\ref{DRE}) and taking the ensemble average of the equation, we obtain an equation for the average concentration, namely
\beq
\label{meanconc}
\frac{d\overline{C}}{dt} = -k \left ( \overline{C}^2  + \overline{C^\prime_AC^\prime_B} \right ),
\eeq
where the last term, $\overline{C^\prime_A(t,x) C^\prime_B (t,x)}$, is the cross-covariance between concentration deviations. For a highly-mixed system, this term is negligible compared to the other terms in the equation, but this is no longer the case when species segregation occurs. During segregation, $A$ and $B$ typically develop anti-correlation properties, so this term becomes negative, slowing the mean reaction rate.
Ultimately, we want to solve for the average concentration $\overline{C}$, but to do this we must first derive an explicit expression for the cross-covariance.
With this in mind, we subtract the equation for the mean concentration (\ref{meanconc}) from (\ref{DRE}) to obtain an evolution equation for the concentration deviation
\beq
\label{DRE_fluctuation}
\frac{\partial C^\prime_i}{\partial t} - D\Delta C^\prime_i = -k\overline{C} (C^\prime_A + C^\prime_B) - kC^\prime_AC^\prime_B + k\overline{C^\prime_AC^\prime_B}, \qquad i = A, B.
\eeq

\subsection{Moment equations and initial conditions}
\label{sec:moment_eqn_IC}
Equations for the evolution of the autocovariance $f(t,x,y) =\bar{ C^\prime_i(t,x)C^\prime_i(t,y)}$, $i=A,B$, and the cross-covariance $g(t,x,y) =\bar{ C^\prime_A(t,x)C^\prime_B(t,y)}$ can be obtained from (\ref{DRE_fluctuation}). Details of the derivation are given in \textit{Paster et al.} \cite{Paster_JCP}. The resulting moment equations are
\beq
\bal
\frac{\partial f}{\partial t} - 2D\Delta f&= -2k \left [\overline{C}(t) (f+ g) + h_1(t,x,y) \right ],\\
\frac{\partial g}{\partial t} - 2D\Delta g&= -2k \left [\overline{C}(t) (f+g) + h_1(t,x,y) \right ],
\eal
\eeq
where
$h_1(t,x,y) = \overline{C^\prime_A(t,x)C^\prime_B(t,x)C^\prime_A(t,y)}$
is a third order moment. 
Taking the sum and difference of the moment equations yields two independent equations for $f+ g$ and
$f -g$, respectively.  These are
\beq
\label{system}
\begin{gathered}
\frac{\partial}{\partial t} (f+g) - 2D\Delta (f+g)= -4k \left [\overline{C}(t) (f+ g) + h_1(t,x,y) \right ],\\
\frac{\partial}{\partial t} (f-g) - 2D\Delta (f-g)= 0.
\end{gathered}
\eeq
Assuming that $\vert \overline{C}(t)(f + g)\vert \gg \vert h_1 \vert$, we can neglect the
third order moment and arrive at an approximate solution to the system of PDEs. We shall return to this assumption in Section \ref{sec:results} by using
results from numerical simulations to examine its validity.

Various initial conditions could be assumed for $f$ and $g$. 
Since we are ultimately interested in showing the correlation between numerical particle models and the analytical solutions of the system, we will focus on synthetic initial conditions that match those imposed in the numerical simulations, given Dirac delta or Gaussian kernel representations (see \ref{appa}). For simplicity and to match numerical results, we assume that the initial cross-covariance is identically zero for each of the different systems (denoted by subscript), namely
\bln$$g_\delta(0,x,y) = g_G(0,x,y) \equiv 0,$$\eln
whereas the autocovariance for the Dirac delta particles is
\beq
\label{AC1}
f_\delta(0,x,y) =C_0m_\delta\left[\delta(x-y)-\frac{1}{\Omega}\right],
\eeq
and for the Gaussian particles is
\beq
\label{AC2}
f_G(0,x,y) = C_0m_G\left[\frac{1}{(4\pi\ell_G^2)^{d/2}}e^{-\frac{\vert x-y\vert^2}{4\ell_G^2}}-\frac{1}{\Omega}\right],
\eeq
where $m_\delta$ and $m_G$ are the mass of single Dirac or Gaussian particle, respectively, $\ell_G$ is the half-width of a Gaussian kernel, and $\Omega$ is the size of a $d$-dimensional domain in a corresponding particle tracking simulation (\ref{appa}).

For these initial conditions and neglecting the effects of higher-order moments like $h_1$, the solution of (\ref{system}) gives rise to a time-dependent cross-correlation function for particle types $p = \delta,G$ (see \textit{Paster et al.} \cite{Paster_JCP}, Appendices A and B for details regarding this derivation):
\beq
\label{g}
g_p(t) = \frac{1}{2(8\pi Dt)^{d/2}}  \int f_p(0,z,x) e^{-\frac{\vert x - z \vert^2}{8Dt}}dz \left [-1 + \exp \left ( -4k\int_0^t \overline{C}_p(\tau) \ d\tau \right ) \right ].
\eeq
\aside{We note that if $f_p(0,z,x)$ is a function of $\vert x-z\vert$ only (as in \eqref{AC1} and \eqref{AC2}), then $g_p$ is independent of $x$, and thus, we write $g_p(t)$ instead of $g_p(t,x,x)$.} Finally, \eqref{g} can be substituted into (\ref{meanconc}) to obtain a closed-form, integro-differential equation for the mean concentration
\beq
\label{cbar_ode}
\frac{d\overline{C}_p}{dt} = -k \overline{C}_p^2 - \frac{k}{2(8\pi Dt)^{d/2}}  \int f_p(0,z,x) e^{-\frac{\vert x - z \vert^2}{8Dt}}dz \left [-1 + \exp \left ( -4k\int_0^t \overline{C}_p(\tau) \ d\tau \right ) \right ],
\eeq
subject to the initial condition $\overline{C}_p(0) = C_0$, with $p = \delta, G$. An explicit representation for the solution of this equation cannot be obtained, but the equation can be solved numerically. 

\section{Particle tracking model}\label{sec:particle}

{\em Benson and Meerschaert} \cite{Benson_react} developed their reactive particle tracking (RPT) algorithm in order to stochastically simulate the DRE using a large number of point-particles (numerically-represented by Dirac deltas). {\em Paster et al.} \cite{Paster_JCP} showed that the initial number of particles necessary for accurate simulation is based on the initial auto- and cross-covariance structure of the system of interest. These authors also show that for systems with either high concentration variance or high correlation of concentration fluctuations, RPT simulations require a relatively small number of particles, while in short-range, low variance systems, the number of particles required to capture this behavior tends toward infinity, making simulations computationally infeasible.

The RPT algorithm uses a particle-killing approach, wherein particles are eliminated from the system after reaction. However, {\em Bolster et al.} \cite{Bolster_mass} discuss a particle-number conserving, mass-based version of RPT (mRPT), in which all particles are given an initial mass that is reduced based on its probability of reaction with every particle of opposite species in the domain. As a result, mRPT is able to depict concentrations of reactant with a higher level of resolution but still requires a uniquely-determined number of point-particles, informed by the initial conditions.

This paper outlines a kernel-based RPT (kRPT) method, which allows the user to reduce the number of particles in a given simulation while still achieving similar results. This is done by ``smearing" the mass of the point particles in mRPT using a kernel representation (Gaussian kernels are used for computational convenience) and thereby smoothing the concentration profile that was coarsened by the decrease in particle number. However, of particular concern is the level of smoothing desired (determined by the half-width of the Gaussian kernel), as an inappropriate choice of kernel size will lead to over- or underprediction of the rate of reaction. The method for determining kernel size, based on the desired number of particles, is discussed in detail within Section \ref{var_calc}.

\subsection{Algorithm details}

The kRPT algorithm initializes the domain by assigning initial particle positions according to draws from a uniform distribution. Diffusion and reaction are simulated separately by operator splitting within each time step, and, without loss of generality, we will assume reactions are calculated first. Reactions are calculated according to the mRPT algorithm of \textit{Bolster et al.} \cite{Bolster_mass} (with the significant difference lying in the increased co-location probability for Gaussian particles, due to the kernel's increased spatial spread). Reactions are performed \aside{for each time step} by reduction in mass and determined probabilistically, based on separation distance, \aside{according to the following algorithm (without loss of generality, as to the choice of inner and outer loop)}
\aside{
\beq
\bal
&\text{\textbf{for}}\ j=1:N_A\\
&\qquad \text{\textbf{for}}\ l=1:N_B\\
&\qquad\qquad \Delta m_{j,l}=k\Delta t\ m_{A,j}m_{B,l}\ v(x_{A,j}-x_{B,l};\Delta t)\\
&\qquad\qquad m_{A,j}=m_{A,j}-\Delta m_{j,l}\\
&\qquad\qquad m_{B,l}=m_{B,l}-\Delta m_{j,l}\\
&\qquad \text{\textbf{end}}\\
    &\text{\textbf{end}}
\eal,
\eeq
such that
\beq
\bal
m_{A,j}(t+\Delta t)&=m_{A,j}(t)-\sum_{l=1}^{N_B}\Delta m_{j,l}(t),\\
m_{B,l}(t+\Delta t)&=m_{B,l}(t)-\sum_{j=1}^{N_A}\Delta m_{j,l}(t),
\eal
\eeq
}
where $m_{A,j}$, $m_{B,l}$ and $x_{A,j}$, $x_{B,l}$ are the masses and positions of the $j^{\text{th}}$ and $l^{\text{th}}$ particle of \emph{A} and \emph{B} reactants, respectively, $N_i$ is the number of particles ($i=A,B$), $D$ is assumed to be the same for \emph{A} and \emph{B}, $\Delta t$ is the constant length of a time step, and $v(s;\Delta t)$ represents the co-location probability of an individual \emph{A} and \emph{B} particle pair, based on their separation distance $s$, given the parameter $\Delta t$. This colocation probability is calculated as follows
\beq
\label{v_s}
v(s;\Delta t)=\frac{1}{\left[4\pi\left(\ell_G^2+2D\Delta t\right)\right]^{d/2}}\exp\left[{-\frac{\vert s\vert^2}{4\left(\ell_G^2+2D\Delta t\right)}}\right],
\eeq
where $\ell_G$ is the half-width of the Gaussian kernel and $2D\Delta t$ is the variance of the diffusion kernel, based on Brownian motion. The positions $x_{A,j}$ and $x_{B,l}$ evolve according to the stochastic Langevin equation
\beq
x_n(t+\Delta t)=x_n(t)+\xi_n\sqrt{2D\Delta t},
\eeq
where $x_n$ is a vector of particle locations, and $\xi_n$ is a vector of independent standard Normal variables, as derived from Fick's Law.

\subsubsection{\aside{Convergence of algorithm to DRE}}
\label{DRE_convergence}

\aside{In order to ensure that the kRPT algorithm is valid for solving diffusion-reaction problems, it would be desirable to show that it converges to the underlying DRE as $\Delta t\to0$. \textit{Paster et al.} \cite{Paster_WRR} demonstrate this convergence for the particle-killing method (RPT) that is the predecessor of mRPT. Additionally, \textit{Bolster et al.} \cite{Bolster_mass} show a similar proof for the mRPT algorithm that is the Dirac particle analog of kRPT. Considering the proof of \textit{Bolster et al.}, it is clear that the proof in the kRPT case follows similarly for any non-Dirac kernel that is (a) a valid density (i.e., the kernel integrates to unity), and (b) both radially-symmetric and even with respect to the origin in the radial variable. As a result, using the proposed kernels with colocation probability $v(s)$, given by \eqref{v_s}, would result in the kRPT algorithm converging to the DRE as the defined time step, $\Delta t$, tends to zero.
}

\subsection{Determination of Gaussian kernel variance}
\label{var_calc}

As defined in \ref{appa}, Gaussian particles will be represented in space by the following kernel
\beq
\label{gauss_kernel}
\phi_G(x-y)=\frac{1}{(2\pi\ell_G^2)^{d/2}}e^{-\frac{\vert x-y\vert^2}{2\ell_G^2}},
\eeq
where $\phi_G$ is a radially-symmetric Gaussian (i.e., half-width is equal in each direction).

Initial concentrations in the system are fixed, and thus the total mass is as well with $C_0={m_{total}}/{\Omega}$, and $m_{total}=m_pN_p$, where $m_p$ is the mass of a single particle of $N_p$ total particles ($p=\delta,G$). This means that once the desired number of particles in the Gaussian simulation ($N_G$) is chosen, then $m_G$ is determined. In order to solve for $\ell_G$ in \eqref{gauss_kernel}, the Gaussian kernel's half-width, we turn to the underlying equation for $\bar C(t)$ given in \eqref{cbar_ode}
\beq
\label{cbar_ode_simple}
\frac{d \bar C_p}{dt}=-k\left(\bar C_p^2+g_p(t)\right).
\eeq

Equation \eqref{cbar_ode_simple} is subject to the following initial conditions, which should match the starting conditions of the particle methods. \aside{Thus, we define the initial conditions, $f_\delta(0,x,y)$ and $f_G(0,x,y)$, to be equal to $\hat f_\delta(0,x,y)$ and $\hat f_G(0,x,y)$, the autocovariance functions for the particle simulations, as derived in \ref{appa}:}

\beq
\label{IC_cont}
\bal
f_\delta(0,x,y):&=\sigma_\delta^2\ell_\delta^d\left[\delta(x-y)-\frac{1}{\Omega}\right]\\
&=\frac{C_0^2\Omega}{N_\delta}\left[\delta(x-y)-\frac{1}{\Omega}\right]\\
&=C_0m_\delta\left[\delta(x-y)-\frac{1}{\Omega}\right]=\hat f_\delta(0,x,y),\\
f_G(0,x,y):&=C_0m_G\left[\frac{1}{(4\pi\ell_G^2)^{d/2}}e^{-\frac{\vert x-y\vert^2}{4\ell_G^2}}-\frac{1}{\Omega}\right]=\hat f_G(0,x,y),\\
g_{p}(0,x,y):&\equiv0,\quad p=\delta,G,
\eal
\eeq
because as established by {\em Paster et al.} \cite{Paster_JCP}, $\sigma_\delta^2\ell_\delta^d={C_0^2\Omega}/{N_\delta}$. Na\"ively, one might attempt to match a Gaussian simulation to a corresponding Dirac simulation by matching properties of the autocovariance functions, $f_\delta$ and $f_G$ in \eqref{IC_cont}. However, this would imply $m_G=m_\delta$ and $\ell_G=0$, turning the Gaussian into a Dirac delta, and making the result trivial. This implies a more nuanced approach is necessary.

We will show in Section \ref{sec:properror} that the error in $\bar C_G(t)$ for differing kernel choices is almost entirely determined by the quantity $\epsilon:=\vert g_\delta(t)-g_G(t)\vert$ \aside{(see Equation \eqref{Cdiff}).} With the goal of minimizing this quantity, we substitute the initial autocovariance choices from \eqref{IC_cont} into \eqref{g} and arrive at the following equations for the cross-covariance functions
\beq
\bal
\label{g_delta_g}
g_\delta(t)&=\frac{C_0m_\delta\left(e^{-4 k\int_0^t \bar C_\delta(\tau)d\tau}-1\right)}{2}\left[\frac{1}{(8\pi Dt)^{d/2}}-\frac{1}{\Omega}\right],\\
g_G(t)&=\frac{C_0m_G\left(e^{-4k\int_0^t\bar C_G(\tau)d\tau}-1\right)}{2}\left[\frac{1}{\left[4\pi(\ell_G^2+2Dt)\right]^{d/2}}-\frac{1}{\Omega}\right].
\eal
\eeq

Using fixed parameters of the simulation, and our choice of $N_G$, we may choose $\ell_G$ such that $\epsilon$ is minimized. Here, we propose two methods for minimizing $\epsilon$.

\subsubsection{Match at specific time}
\label{sec:tstar}

The most straightforward way of minimizing $\epsilon$ is choosing a particular time, $t^*$, at which the user would like the simulations to ``match," so that $g_G(t^*)=g_\delta(t^*)$, which implies $\epsilon(t^*)=0$. We assume here that $\vert \bar C_\delta-\bar C_G\vert\approx0$ on the interval $[0,t^*]$, a claim supported by Equation \eqref{Cdiffend}. This allows us to cancel the exponential terms in $g_\delta$ and $g_G$ and solve for $\ell_G$ using only the parameters of the system. Thus, we may attain a value for $\ell_G(t^*)$ by setting $g_\delta(t^*)=g_G(t^*)$ in \eqref{g_delta_g} and solving, which gives
\beq
\label{lg_of_tstar}
\ell_G(t^*)=\sqrt{\frac{1}{4\pi}\left[\frac{N_G}{N_\delta}\left[(8\pi D t^*)^{-d/2}-\frac{1}{\Omega}\right]+\frac{1}{\Omega}\right]^{-2/d}-2Dt^*}.
\eeq

Here, we note that the form of \eqref{lg_of_tstar} implies that, for given parameters, there is some maximum value of $t^*$, which we will call $\tau^*$, such that
\beq
\label{tau_star}
\tau^*:=\max\left\{t^*:\left[\frac{N_G}{N_\delta}\left[(8\pi D t^*)^{-d/2}-\frac{1}{\Omega}\right]+\frac{1}{\Omega}\right]^{-2/d}-8\pi Dt^*\geq0\right\},
\eeq
and after which, values for $\ell_G$ will become imaginary, \aside{and hence, non-physical.} However, the value of $\tau^*$ is highly sensitive to the parameters $D$ and $\Omega$, simply implying that certain choices of $t^*$ may not be appropriate for systems with particular Damk\"ohler numbers. Note, also, that \aside{the left side of the inequality in \eqref{tau_star}} is equal to zero only when $N_G=N_\delta$ and greater than zero when $N_G<N_\delta$.

This method of minimization will produce a result in which Dirac- and Gaussian-based moment equation solutions will match at or near $t^*$ and allows the user to decide when the differing simulations will display similar behavior.

\subsubsection{Least squares minimization}
\label{sec:least_squares}

Alternatively, the user may define a set of times, $T^*=\{t_0,t_1,\dots,t_n\}$
, over which the quantity
\beq
\epsilon(\ell_G):=\left(\sum_{k=1}^n\vert \bar C_\delta(t_k)-\bar C_G(t_k;\ell_G)\vert^2\right)^{1/2},
\eeq
may be minimized to determine
\beq
\ell_G^* := \mathop{\mathrm{argmin}}_{\ell_G\in\left(0,\frac{1}{c}\Omega\right]}\quad\epsilon(\ell_G),
\eeq
where the upper-bound for the value of $\ell_G$ is some fraction of the domain, $\Omega$.
This will produce a single corresponding value for $\ell_G = \ell_G^*$, which will minimize the overall error between Dirac- and Gaussian-based moment equation solutions, subject to the chosen set of times $T^*$.

\section{Propagation of error}
\label{sec:properror}

In considering the differing initial autocovariances arising from Dirac versus Gaussian kernels, \aside{as shown in \eqref{IC_cont},}
one arrives at two different solutions of the diffusion-reaction equation for \aside{$g(t)$.}
In either case, the form of $g$ is the same, namely
\beq
\label{g2}
g_p(t) = \psi_p(t) \left [-1 + \exp \left ( -4k\int_0^t \overline{C}_p(\tau) \ d\tau \right ) \right ],
\eeq
where
\beq
\label{psi}
\psi_p(t) = \frac{1}{2(8\pi Dt)^{d/2}} \int f_p(0,z,x) e^{-\frac{\vert x - z \vert^2}{8Dt}} \ dz,
\eeq
for $p=\delta, G$. We note that $\psi$ does not depend upon $x$ because each of the initial conditions (\ref{AC1})-(\ref{AC2}) is a function only of the difference $x - y$, rather than $x$ and $y$ independently, and the integral is translation invariant.
In particular, a few brief calculations show that
\beq
\label{psid}
\psi_\delta(t) = \frac{1}{2}C_0 m_\delta \left [\frac{1}{(8\pi Dt)^{d/2}} - \frac{1}{\Omega} \right ],
\eeq
and
\beq
\label{psiG}
\psi_G(t) = \frac{1}{2} C_0 m_G \left [ \frac{1}{(4\pi(\ell_G^2 + 2Dt))^{d/2}} - \frac{1}{\Omega} \right ].
\eeq

For each $g_p$, we have a corresponding mean concentration $\overline{C}_p(t)$ defined as the solution of the initial-value problem
\beq
\label{Cbar}
\left.
\begin{array}{l}
\displaystyle \frac{d\overline{C}_p}{dt} = -k \left [ \overline{C}_p^2 + g_p(t) \right ]\\
\\
\overline{C}_p(0) = C_0.
\end{array}
\right \}
\eeq
Thus, in order to estimate the error in the resulting mean concentration stemming from the use of a different (and non-Dirac delta) autocorrelation function, we first assume, without loss of generality, that $\overline{C}_G(t) \geq \overline{C}_\delta(t)$.
Therefore, to estimate the difference in these quantities, we subtract (\ref{Cbar}) in the case $p=\delta$ from (\ref{Cbar}) in the case $p=G$ so that
\bln
\begin{eqnarray*}
\frac{d}{dt} \left ( \overline{C}_G - \overline{C}_\delta \right) & = & -k \left ( \overline{C}_G^2 - \overline{C}_\delta^2 \right) - k(g_G - g_\delta)\\
& = & -k(\overline{C}_G + \overline{C}_\delta)(\overline{C}_G - \overline{C}_\delta) + k(g_\delta - g_G).
\end{eqnarray*}
\eln
Using the lower bound for each mean concentration, (\ref{Clower}), in \ref{appb}, this becomes
\bln$$ \frac{d}{dt} \left ( \overline{C}_G - \overline{C}_\delta \right) \leq -\frac{2kC_0}{1 + kC_0 t} (\overline{C}_G - \overline{C}_\delta) + k(g_\delta - g_G).$$\eln
Now, using the integrating factor $p(t) = (1+ kC_0t)^2$, we can rearrange the inequality and integrate using 
$\overline{C}_\delta(0) = \overline{C}_G(0) = C_0$ to find
\bln$$\overline{C}_G(t) - \overline{C}_\delta(t) \leq \frac{k}{(1 + kC_0t)^2} \int_0^t [ g_\delta(\tau) - g_G(\tau) ] (1+ kC_0\tau)^2 \ d \tau.$$\eln
Since these same computations can be repeated on intervals of time in which $\overline{C}_\delta(t) \geq \overline{C}_G(t)$ we combine them with absolute values, resulting in
\beq
\label{Cdiff}
\left \vert \overline{C}_G(t) - \overline{C}_\delta(t) \right \vert \leq \frac{k}{(1 + kC_0t)^2} \int_0^t  \vert g_\delta(\tau) - g_G(\tau) \vert  (1+ kC_0\tau)^2 \ d \tau.
\eeq
Thus, the difference in mean concentration can be estimated in terms of $\vert g_\delta(\tau)-g_G(\tau)\vert$ only.

Next, we need to estimate the difference in the cross-correlation functions caused by the differing initial autocorrelation structure.  This is performed in detail within \ref{appb}, resulting in (\ref{gbound}), namely
\bln$$\vert g_\delta(t) - g_G(t) \vert \leq  \frac{1}{2}C_0^2 \Delta + \frac{C_0^2 \Delta}{2\Omega (1 + kC_0t)^{d/2}} + \frac{\pi dC_0 \ell_G^2m_G}{(1 + kC_0t)^{1+ d/2}} +  \frac{C_0 m_\delta}{2\Omega(1 + kC_0t)^{4}},$$\eln
where
\bln$$ \Delta := \frac{1}{N_G} - \frac{1}{N_\delta} = \frac{1}{N_G} \left ( 1 - \frac{N_G}{N_\delta} \right )$$\eln
can be chosen as small as one wishes.
Inserting this expression for the difference in cross-correlations into (\ref{Cdiff}), we find
\bln
\begin{eqnarray*}
\left \vert \overline{C}_G(t) - \overline{C}_\delta(t) \right \vert & \leq & \frac{k}{(1 + kC_0t)^2} \int_0^t \left [(1+ kC_0\tau)^2 \frac{1}{2}C_0^2 \Delta + \frac{C_0^2 \Delta}{2\Omega (1 + kC_0\tau)^{d/2 - 2}} \right. \\
& \ & \left. + \frac{\pi dC_0 \ell_G^2m_G}{(1 + kC_0\tau)^{-1+ d/2}} +  \frac{C_0 m_\delta}{2\Omega(1 + kC_0\tau)^{2}} \right ] \ d \tau\\
& \leq &  \frac{1}{(1 + kC_0t)^2} \left [ C_1\Delta (1+kC_0 t)^{3} + C_2\Delta (1+kC_0 t)^{3-d/2} \right.\\
& \ & \left. + C_3(1+kC_0 t)^{2-d/2} + C_4 \right ]\\
& = & C_1\Delta (1 + kC_0t) + C_2\Delta (1 + kC_0t)^{1 - d/2} \\
& \ & + C_3(1 + kC_0t)^{- d/2} + C_4(1 + kC_0t)^{-2},
\end{eqnarray*}
\eln
where $C_1, C_2, C_3$, and $C_4$ are constants defined by the terms within the first and second lines above.
For $t$ large, the second, third, and fourth terms are negligible compared to the first.  Hence, denoting the stopping time of a simulation by $T$, the difference in the mean concentration computed by the RPT method and the kRPT implementation satisfies
\beq
\label{Cdiffend}
\left \vert \overline{C}_G(t) - \overline{C}_\delta(t) \right \vert \leq C_1\Delta (1 + kC_0 T).
\eeq
Of course, choosing $\Delta$ sufficiently small will allow for as great an order of accuracy as one wishes.
In particular, if one chooses $\Delta = 0$, which occurs when $N_G = N_\delta$, the error instead satisfies
\beq
\label{Cdiffend2}
\left \vert \overline{C}_G(t) - \overline{C}_\delta(t) \right \vert \leq C_3(1 + kC_0T)^{- d/2},
\eeq
which actually decreases with $T$, and this error is incurred due only to the difference in initial covariance structure.

\section{Results}\label{sec:results}

\begin{figure}[t]
\centering
\includegraphics[width=\textwidth]{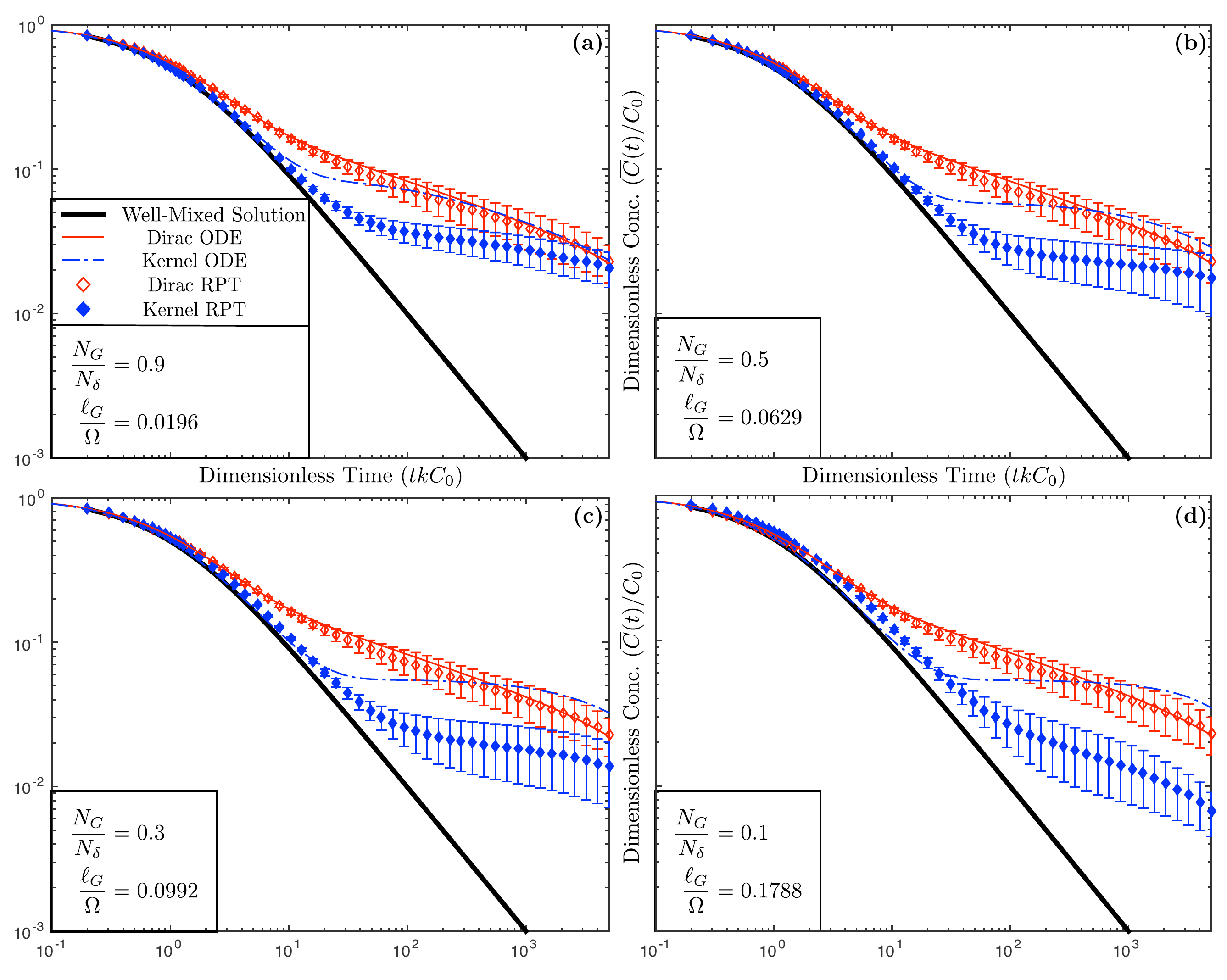}
\caption{Plots of particle tracking simulations (symbols) and moment equation solutions (curves), using specific time matching for $t^*=100$ ($t^*kC_0=500$), for a range of $N_G/N_\delta$ values.}
\label{fig1}
\end{figure}

\begin{figure}[t]
\centering
\includegraphics[width=\textwidth]{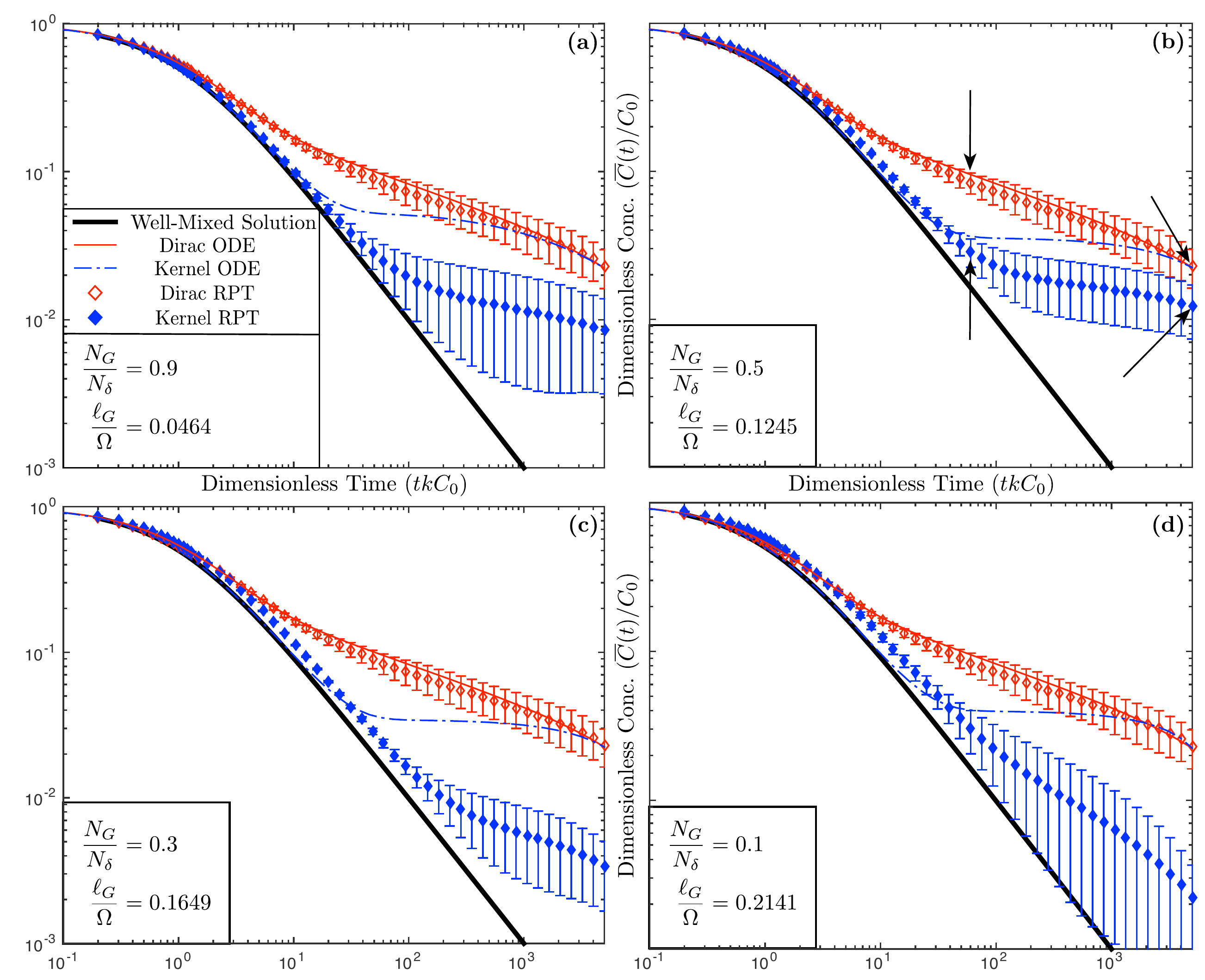}
\caption{Plots of particle tracking simulations (symbols) and moment equation solutions (curves), using specific time matching for $t^*=1000$ ($t^*kC_0=5000$), for a range of $N_G/N_\delta$ values. Arrows in (b) correspond to the plots shown in Figure \ref{fig7}(a)-(d).}
\label{fig2}
\end{figure}

\begin{figure}[t]
\centering
\includegraphics[width=\textwidth]{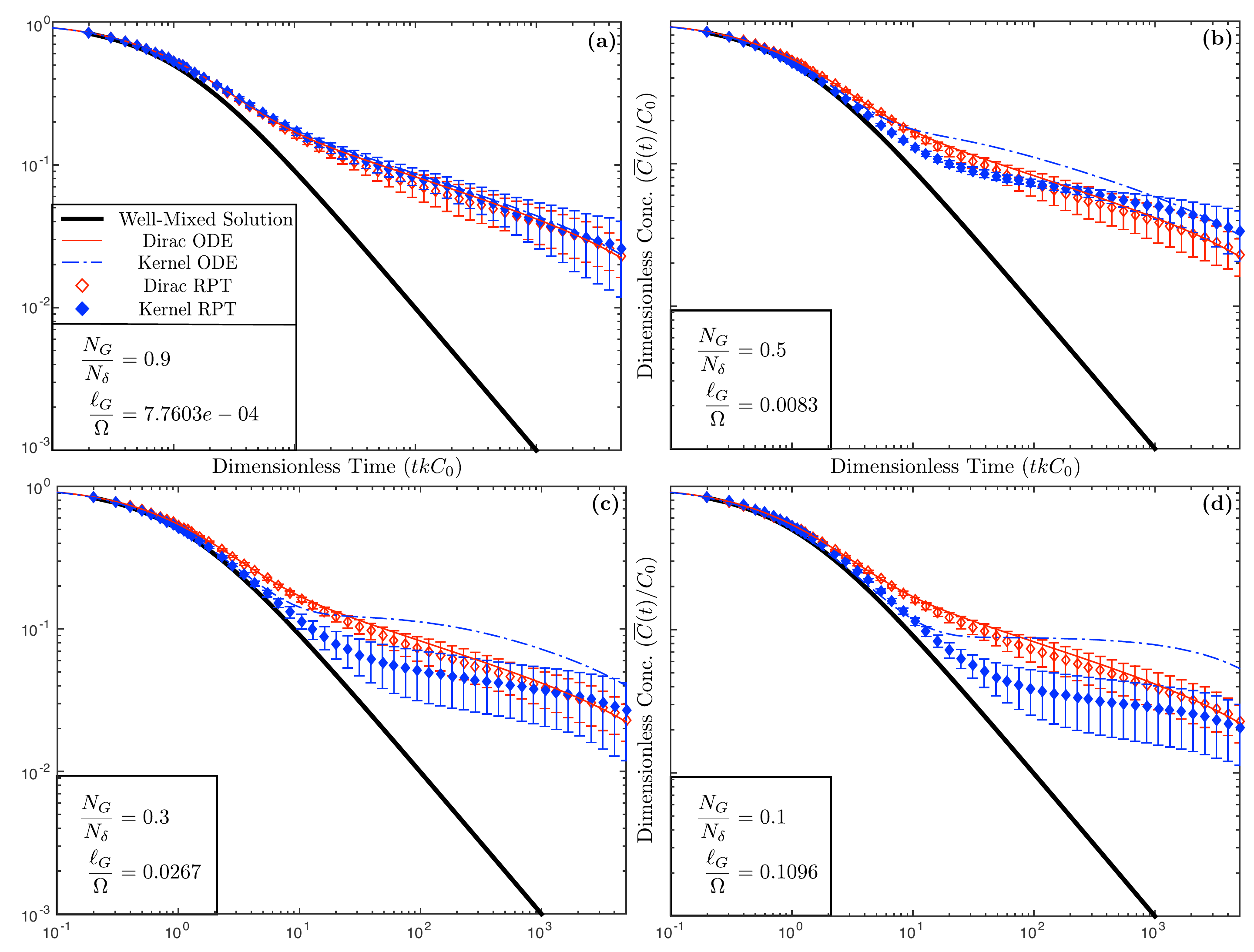}
\caption{Plots of particle tracking simulations (symbols) and moment equation solutions (curves), using least squares matching for a range of $N_G/N_\delta$ values.}
\label{fig3}
\end{figure}

\begin{figure}[t]
\centering
\includegraphics[width=\textwidth]{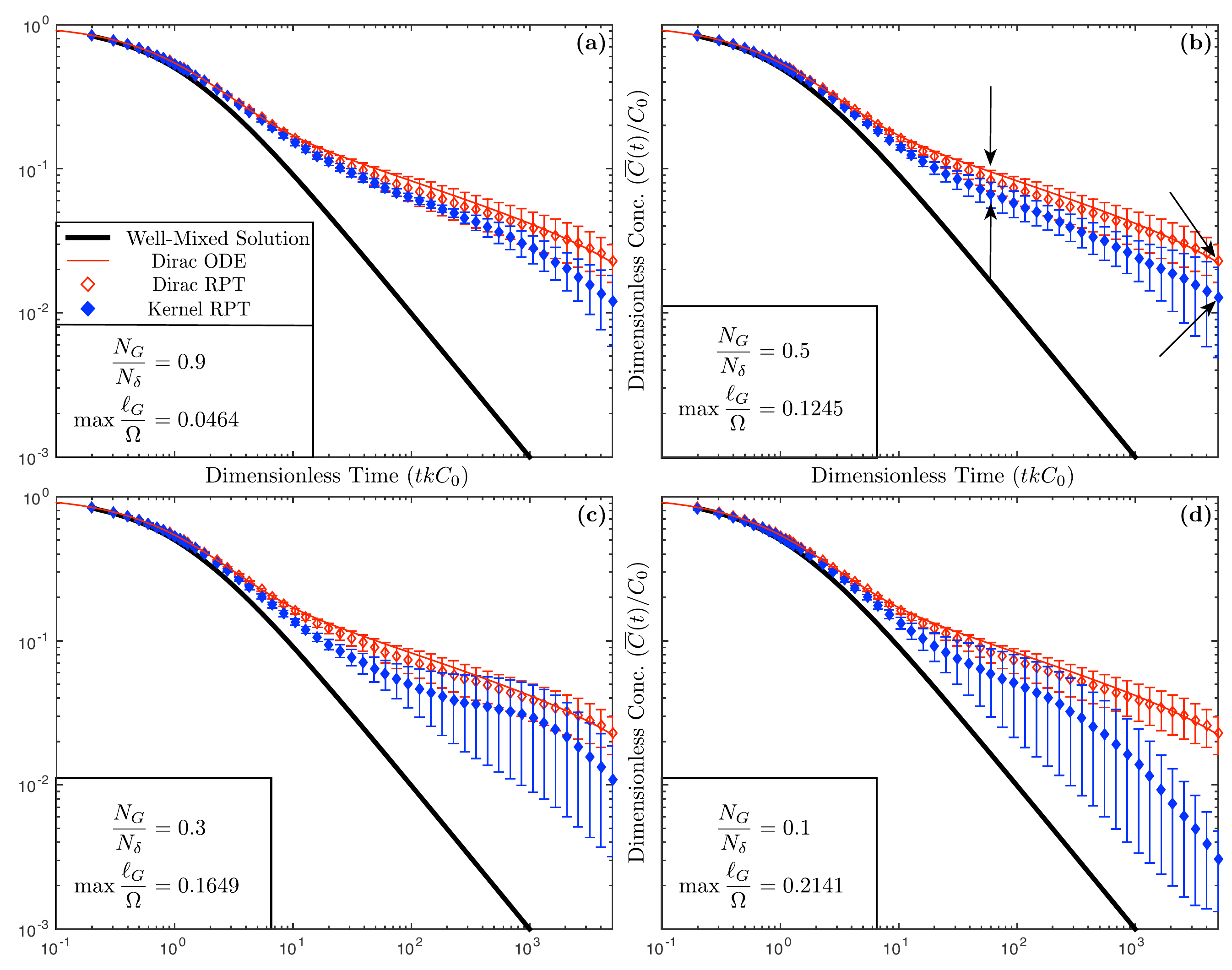}
\caption{Plots of particle tracking simulations (symbols) and moment equation solutions (curves), using variable $\ell_G$ matching for a range of $N_G/N_\delta$ values. Arrows in (b) correspond to the plots shown in Figure \ref{fig7}(a)-(b) and (e)-(f).}
\label{fig4}
\end{figure}

\begin{figure}[t]
\centering
\includegraphics[width=\textwidth]{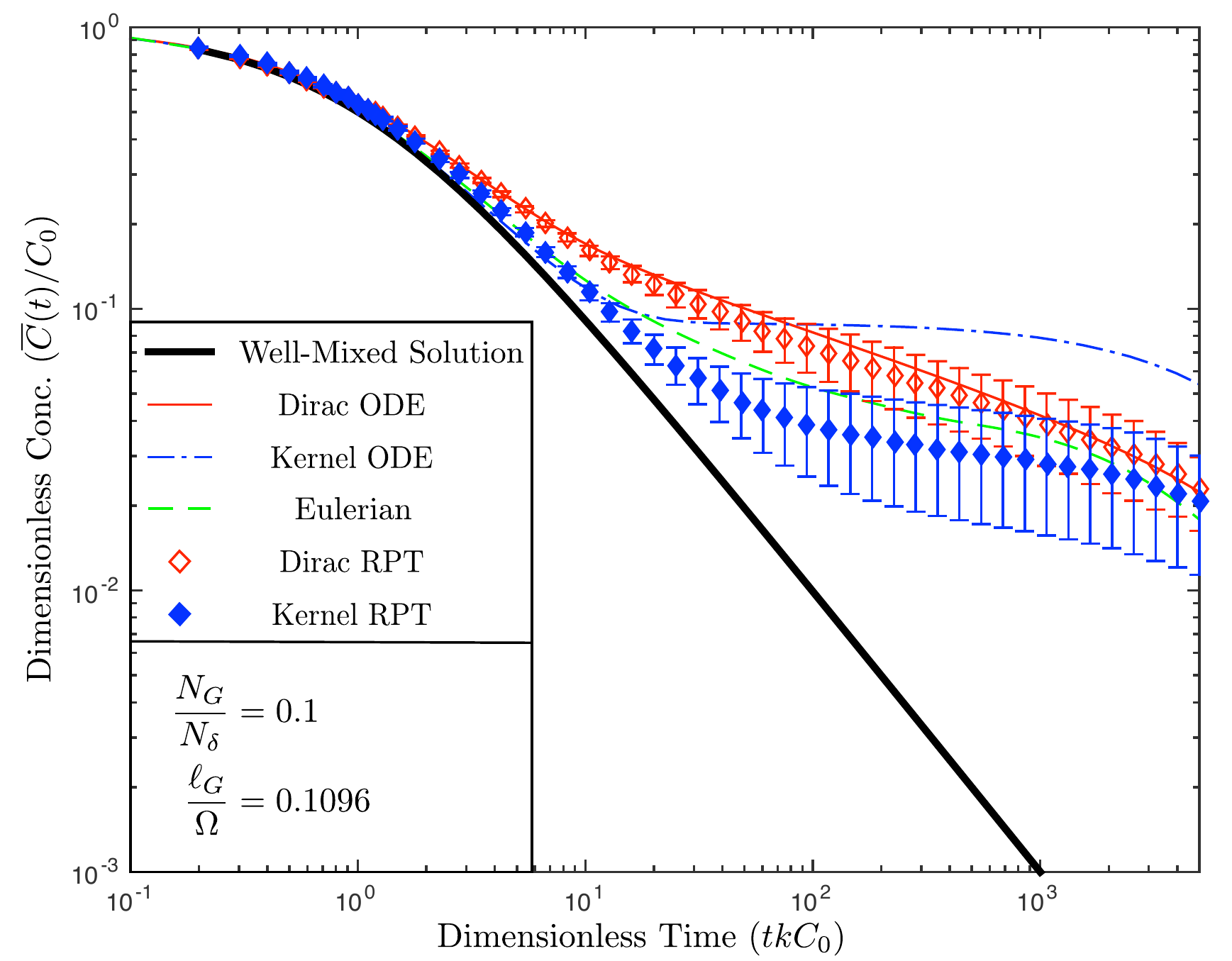}
\caption{\aside{Plots of particle tracking simulations and moment equation solutions as compared to an Eulerian finite difference solution, using the least squares matching criteria for $N_G/N_\delta=0.1$ (corresponding to Figure \ref{fig3}(d)).}}
\label{eulerian}
\end{figure}

\begin{figure}[t]
\centering
\includegraphics[width=\textwidth]{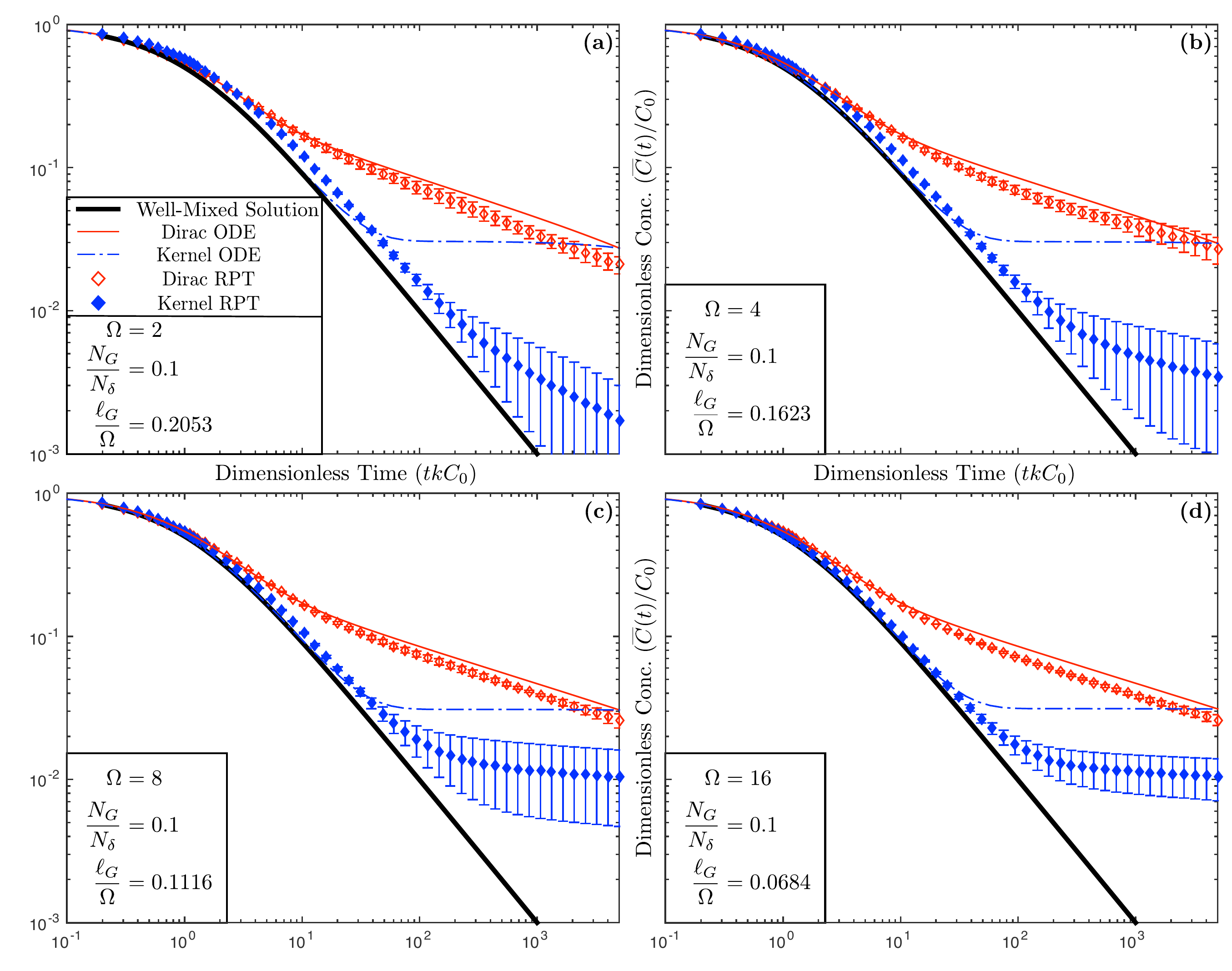}
\caption{Plots of particle tracking simulations (symbols) and moment equation solutions (curves), using specific time matching with $t^*=1000$ ($t^*kC_0=5000$), $N_G/N_\delta=0.1$ for a range of domain sizes.}
\label{fig5}
\end{figure}

\begin{figure}[t]
\centering
\includegraphics[width=\textwidth]{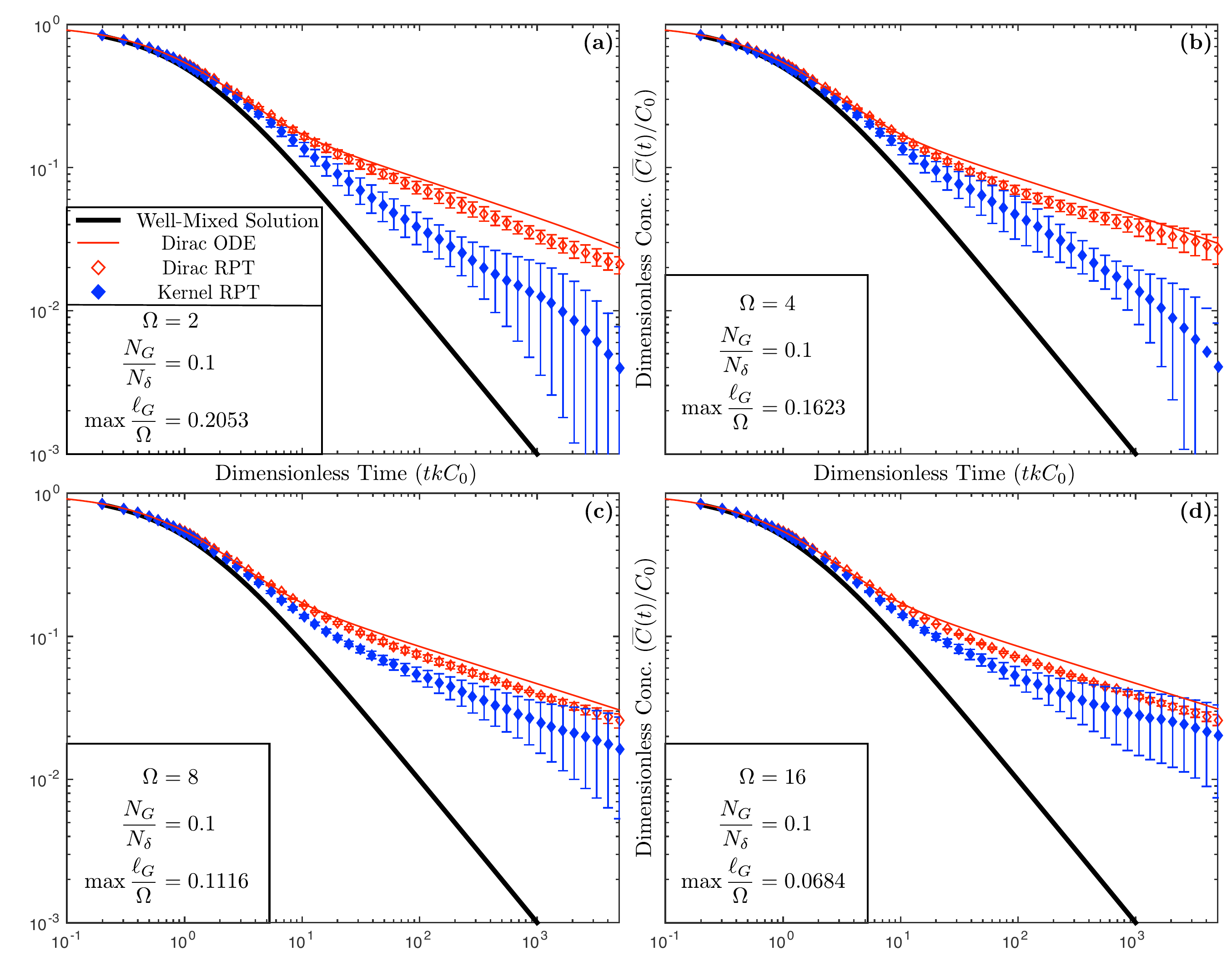}
\caption{Plots of numerical simulations (symbols) and moment equation solutions (curves), using variable $\ell_G$ matching with $N_G/N_\delta=0.1$ for a range of domain sizes.}
\label{fig6}
\end{figure}

\begin{figure}[t]
\centering
\includegraphics[width=.86\textwidth]{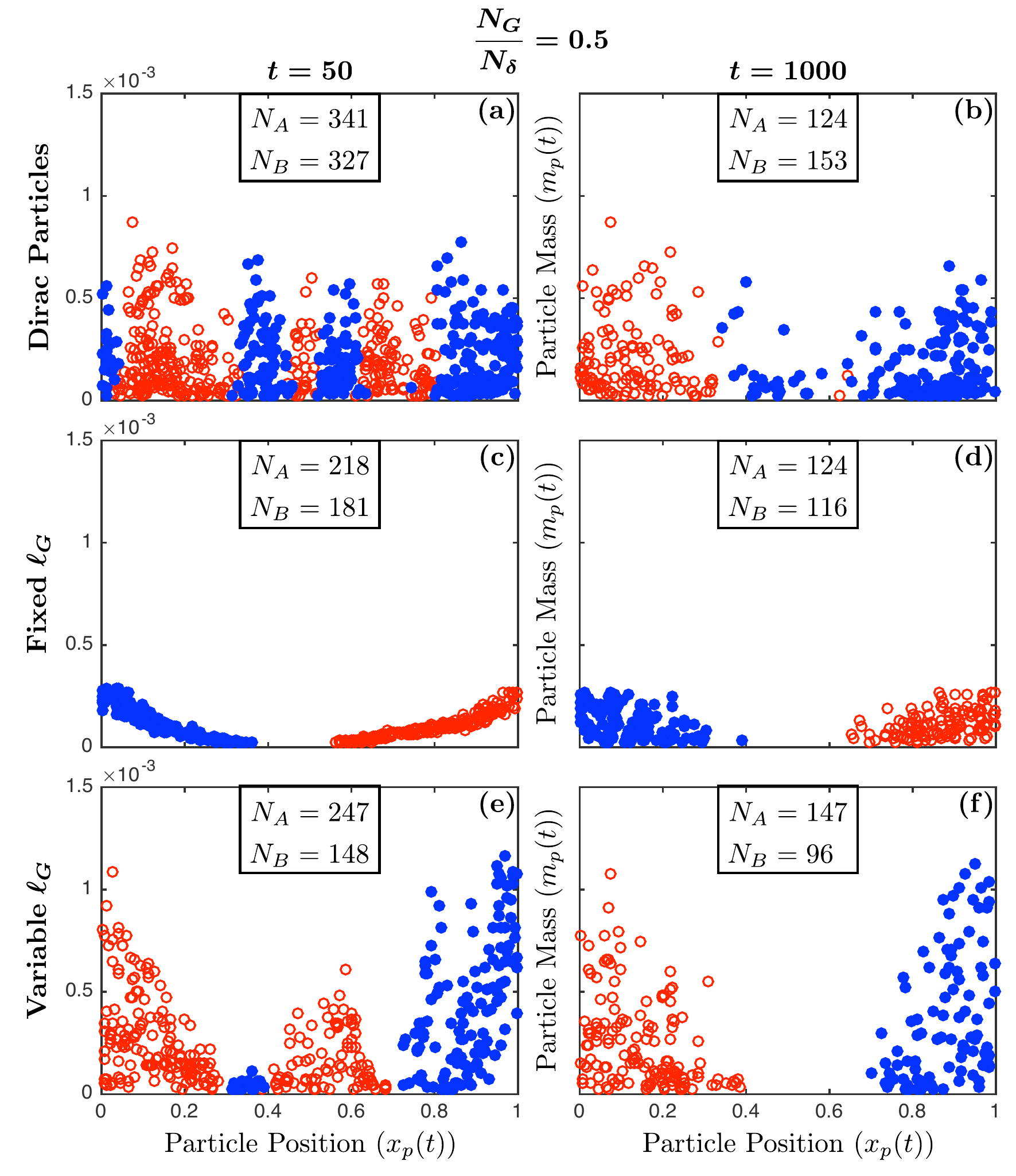}
\caption{Plots of particle position/mass demonstrating moats phenomenon for Gaussian kernels (bottom two rows). Plots show $N_G/N_\delta=0.5$ using: Dirac kernels, specific time matching (fixed $\ell_G$) with $t^*=1000$ ($t^*kC_0=5000$), and variable $\ell_G$ matching.}
\label{fig7}
\end{figure}

Numerical simulations were conducted in MATLAB, using a MacBook Pro with a 2.9 GHz Intel Core i5 processor and 8 GB of RAM. All simulations were conducted for $d=1$ spatial dimensions, and six-realization ensembles were employed in order to average the stochastic variability of particle tracking simulations. In all plots showing domain-averaged concentration versus time, the well-mixed analytical solution \eqref{wellmixed} is shown as a solid black line for reference, moment equation solutions \eqref{cbar_ode} are shown as solid curves, and particle tracking solutions are shown as diamond-shaped scatter plots with error bars corresponding to $\pm1$ standard deviation among realizations in the ensemble. As well, red plots (solid curve, empty diamonds) represent Dirac kernel simulations, and blue plots (dashed curve, filled diamonds) represent Gaussian kernel simulations. The domain size for all initial cases is $\Omega=1$ and is expanded when considering domain effects.

In this section, we characterize the simulations by their particle Damk\"ohler number, as defined in \textit{Bolster et al.} \cite{Bolster_mass}, $\widehat{Da}:=kC_0\Delta x^2/D$, where $\Delta x$, the average inter-particle spacing ($\Delta x=\Omega/N_p$, $p=\delta,G$), is used as the characteristic length-scale of diffusive mixing. Note that this characterization of the system is derived using Dirac particles, and, for all numerical simulations in this section, $\widehat{Da}=0.5$ will describe the Dirac particle base cases ($D=1.0\times10^{-5}$, $k=5.0$, $C_0=1.0$).

Throughout this section, we will show the results of numerical simulations for $N_G=\{900,500,300,100\}$ and compare them to the base case of $N_\delta=1000$. Solution matching criteria used will be: specific time matching for $t^*=\{100,1000\}$, least squares matching, and a heuristically-motivated variable kernel half-width matching. Additionally, the various influences of the ratio of kernel half-width to domain size ($\ell_G/\Omega$) will be examined for several sub-cases.

\subsection{Specific time matching}

As the simplest and most straightforward approach to matching Dirac and Gaussian particle tracking solutions, we first use the specific time matching criteria put forth in Section \ref{sec:tstar}. Simulations were run for $t^*=\{100,1000\}$. As depicted in Figures \ref{fig1}(a) and \ref{fig2}(a), we predictably see very close agreement between the particle tracking and moment equation solutions for $N_G/N_\delta=0.9$. As well, in both cases, the Dirac and Gaussian particle tracking solutions become ``close" near $t^*$ (the error in dimensionless concentration between the Dirac and particle tracking solutions at $t^*$ is $\mathcal{O}(1/100)$, or approximately the initial mass of 10 Gaussian particles), and remain close thereafter in the $t^*=100$ case. However, one does notice the overprediction of concentration by the moment equation solutions (blue curves), versus the particle tracking solutions for the Gaussian kernels. In the cases of $N_G/N_\delta=0.9$ (Figures \ref{fig1}(a) and \ref{fig2}(a)), we attribute this to the neglect of third-order moments in the moment equation paradigm, as mentioned in Section \ref{sec:moment_eqn_IC}. These third-order moments appear to become more significant as $N_G$ decreases (and $\ell_G$ increases), causing an increasing discrepancy between moment equation and particle tracking solutions as $N_G/N_\delta$ decreases. For a discussion of the influence of third-order moments in a non-particle-number preserving, Dirac kernel case, see \textit{Paster et al.} \cite{Paster_JCP}, though we note here that the influence of these moments appears to be significantly greater when using Gaussian kernels\aside{, due to the difference in initial conditions.}

Considering Figures \ref{fig1}(b)-(d) and \ref{fig2}(b)-(d), we notice a few other behaviors of interest, the effects of which increase with the ratio $\ell_G/\Omega$ (which, in turn, increases with $t^*$). First, we see an increasing separation at early times between the Gaussian particle tracking solution and the Dirac particle solution as $\ell_G/\Omega$ increases. This phenomenon propagates through time, resulting in an overprediction of concentration at late times. Additionally, two late-time behaviors become apparent as $\ell_G/\Omega$ increases. First, concentrations are increasingly underpredicted (which, to some extent, offsets the aforementioned overpredicting behavior, as can be seen in Figure \ref{fig2}(b), in that it shows a slightly closer match at $t^*=1000$ than does Figure \ref{fig2}(a)). Second, as $\ell_G/\Omega$ increases, there emerges a behavior wherein the slope of the Gaussian particle simulation becomes increasingly negative, and is seen most clearly in Figure \ref{fig2}(c) and (d). This behavior is explained in \textit{Bolster et al.} \cite{Bolster2012} as when the particle island sizes are comparable to the domain size. These domain effects will be discussed in greater detail in Section \ref{sec:domain_effects}.

We notice that, using the specific time matching criteria, our moment equation and particle tracking solutions show different behavior, namely the particle solutions not reliably intersecting near $t^*$. As an improvement, the least squares matching criteria has the potential to eliminate some of the unrecoverable errors that begin long before $t^*$.

\subsection{Least squares matching}

For these simulations, the least squares error was minimized over a set of times, $T^*$, with 100 logarithmically-spaced points between $t=10^{-2}$ and 1000, according to the algorithm in Section \ref{sec:least_squares}. While not quite as simple to implement as the specific time matching approach, we see in Figure \ref{fig3}, the least squares matching criteria shows the ability to match Dirac and Gaussian particle solutions with a high degree of fidelity. For example, the maximum error in dimensionless concentration between the particle tracking solutions in the $N_G/N_\delta=0.1$ case (shown in Figure \ref{fig3}(d)) is 0.0498, and final time error is 0.0023 or approximately the initial mass of five and 0.2 Gaussian particles, respectively.

Additionally, the small $\ell_G$ values predicted by the least squares matching (e.g., $\ell_G=0.1096$ for $N_G/N_\delta=0.1$, as in Figure \ref{fig3}(d), corresponds to $t^*\approx 15$ ) causes the above-mentioned domain effects to be an apparent non-issue for the simulated parameters. However, the effect of neglecting third-order moments is unavoidable, as is evidenced by the increasing overprediction in concentration by the moment equations as $N_G/N_\delta$ decreases. This effect is depicted most clearly in Figure \ref{fig3}(d), where the Gaussian particle solution closely matches both the particle and moment equation Dirac solutions, but the Gaussian moment equation solution differs significantly.

\subsection{Variable kernel half-width matching}
\label{sec:var_lg}

Another study by \textit{Rahbaralam et al.} \cite{Dani_kernel} proposed adjusting kernel size as a function of time in a simulation. We propose to treat $\ell_G$ as a function of time, by substituting the simulation time, $t$, for $t^*$ in \eqref{lg_of_tstar}. Additionally, substituting \eqref{lg_of_tstar} into \eqref{g_delta_g} and using this newly formulated cross-covariance in \eqref{cbar_ode_simple}, we get a complicated equation for $\bar C_G(t)$ that may be solved numerically. The expected result prevails that $\bar C_\delta-\bar C_G\equiv0$, as we are simply creating a condition under which $g_\delta$ and $g_G$ are numerically equal for all chosen time steps. We admit here, and leave as an open question, that we have not generated the analytical framework to properly analyze this reformulated moment equation. Nonetheless, numerical particle tracking simulations were conducted using this method, wherein $\ell_G$ is recalculated at every time step, using actual simulation time ($\ell_G$ is initialized using $t=\Delta t/2$), giving good overall matching to the Dirac solutions, as shown in Figure \ref{fig4}. Figure \ref{fig4}(a) and (b) show very close agreement between the Dirac and Gaussian particle simulations, while Figure \ref{fig4}(c) and (d) appear to show the effects of large $\ell_G$ relative to $\Omega$.

In order to examine the capabilities of this method, simulations were conducted that attempted to decrease the domain effects (the details of which are discussed in Section \ref{sec:domain_effects}). As shown in Figure \ref{fig6}, as the effect of the domain is decreased, the variable half-width kernel method converges very closely to the Dirac particle simulations (e.g., Figure \ref{fig6}(d) has a maximum error in dimensionless concentration of 0.0236, and an error of 0.0058 at final time, or approximately the mass of 2 or 0.5 Gaussian particles, respectively).

\subsection{\aside{Comparison to Eulerian method}}
\label{sec:Eulerian}

\aside{In order to further validate the kRPT model, we compare the numerical results to a classical Eulerian numerical method. In this direction, we follow a similar approach to that of \textit{Bolster et al.} \cite{Bolster_mass}, in that initial concentrations are stochastically perturbed about $C_0$ by drawing from a uniform distribution. The method used will be a semi-implicit finite difference scheme. The linear portion of our method (diffusion) will be backward in time, central in space, and the nonlinear portion (reaction) will be backward in time, formulated as follows}
\aside{
\beq
\bal
\frac{C_i(x,t)-C_i(x,t-\Delta t)}{\Delta t}-&D\frac{C_i(x+\Delta x,t)-2C_i(x,t)+C_i(x-\Delta x,t)}{\Delta x^2}=\\
&\hspace{9em}-kC_A(x,t-\Delta t)C_B(x,t-\Delta t),\\
i=&A,B,
\eal
\eeq}
\aside{where $\Delta x$ is chosen to be $1/N_\delta$ and $\Delta t$ is chosen so as to match the time steps of the particle tracking simulations.}

\aside{We compare the finite difference solution to the results of the least squares matching for $N_G/N_\delta=0.1$ (the least squares plot shown in Figure \ref{fig3}(d)), within Figure \ref{eulerian}, which shows a very close agreement between the solutions generated by both the Dirac and Gaussian particle tracking simulations and the Eulerian finite difference simulation.}


\subsection{Domain effects}
\label{sec:domain_effects}

In each of the three methods to estimate $\ell_G$, we see characteristic behaviors of finite domain effects. These include the early time overprediction, late time underprediction, and the increasingly negative slope at late time (for a further discussion of domain effects, see \textit{Bolster et al.} \cite{Bolster2012}). We choose two cases that seemed to suffer from the most significant domain effects, namely the specific time matching case for $N_G/N_\delta=0.1$, $t^*=1000$ and the variable kernel half-width case for $N_G/N_\delta=0.1$ (shown in Figures \ref{fig2}(d) and \ref{fig4}(d)).

In order to hold the Damk\"ohler number of the systems constant as we increase the domain size from $\Omega=1$ to $\Omega=\{2,4,8,16\}$, we increase the particle number in direct proportion. This allowed us to simulate the same system (holding $\widehat{Da}$ constant) while decreasing the values for $\ell_G/\Omega$. Note that since $\Omega$ appears in \eqref{lg_of_tstar}, this relationship is somewhat complicated, but $\ell_G$ increases approximately proportional to $\sqrt{\Omega}$.

As shown in Figures \ref{fig5} and \ref{fig6}, one can readily see the aforementioned domain effects decrease as $\ell_G/\Omega$ decreases. In Figure \ref{fig5}, as the domain expands, the Gaussian particle solutions begin to converge more closely to the well-mixed analytical and Dirac particle solutions at early time. As well, in both Figures \ref{fig5} and \ref{fig6}, the sharp turndown in slope at late time disappears as the kernels no longer ``feel" effects from the boundary, and the concentration underprediction is reduced. The combination of these results leads to significantly improved agreement in both of the considered cases. While we propose no exact proof for the maximum appropriate ratio $\ell_G/\Omega$ to minimize domain effects, our numerical results suggest that a good rule of thumb is to keep this ratio below 0.12.

\subsection{Slowdown in Gaussian particle reaction speed}
\label{sec:moats}

Examining the results of the above numerical simulations, a particular behavior becomes apparent, in that the Gaussian particle simulations tend to have an accelerated reaction rate (as compared to Dirac simulations) at early time, followed by a period during which reaction slows significantly, and finally, at late time, reaction rate increases again, and the slope of the Gaussian particle solution begins to match that of the Dirac solution. This behavior is demonstrated most clearly in Figure \ref{fig2}(b). 

In order to analyze this behavior, single-realization simulations were run where position and mass were recorded for all particles as they evolve in time, and results are shown in Figure \ref{fig7} (only particles with mass greater than $2\%$ of the original mass of one Dirac particle are depicted). These simulations were conducted to compare the Dirac particle simulation to the Gaussian particle case of $N_G/N_\delta=0.5$ in the two subcases of specific time matching for $t^*=1000$ (``fixed $\ell_G$") and variable kernel half-width matching (shown in Figures \ref{fig2}(b) and \ref{fig4}(b)).

What we observe in the Dirac case is the well-documented and  predictable segregation of species and the formation of ``islands" due to low levels of diffusion, as shown in Figure \ref{fig7}(a) and (b), and discussed in \textit{Toussaint and Wilczek} \cite{Toussaint}. This segregation appears to begin around $t=50$, with the formation of small islands (Figure \ref{fig7}(a)), and the segregation becomes quite pronounced at the final time $t=1000$, when we are left with only two large islands (Figure \ref{fig7}(b)). However, in the fixed $\ell_G$ case, we see strict segregation beginning around $t=50$, with the formation of ``moats" around the islands with a width approximately $2\ell_G$ (Figure \ref{fig7}(c)). These islands and moats persist until late time (Figure \ref{fig7}(d)), when the overall reaction is dictated by the time required for particles to migrate to the island edges.

It was this phenomenon that was the motivation for the variable $\ell_G$ approach discussed in Section \ref{sec:var_lg}, as a kernel that starts small and grows with time should prevent the early formation of moats (and corresponding slowdown of reaction rate), while still allowing for solution matching as time progresses and the kernel grows. We see evidence of exactly this at $t=50$ (Figure \ref{fig7}(e)), where smaller islands have begun to form but no moats are apparent ($\ell_G(t=50)=0.0473$). Then, at $t=1000$ (Figure \ref{fig7}(f)), we see a similar moat size to the fixed $\ell_G$ case, but with greater overall particle masses, explaining the closer fit to the Dirac solution in the variable $\ell_G$ case.

\section{Conclusions}\label{sec:concl}

In this paper, we present a kernel-based particle tracking method for modeling imperfectly-mixed chemical reactions in diffusive media. In particular, Gaussian kernels were chosen for computational convenience, notably the preservation of their Gaussian structure under convolution, although any kernel with similar properties could presumably be used. The primary advantage of kRPT is the ability to significantly reduce the number of particles used, as compared to a corresponding Dirac kernel-based mRPT simulation, while introducing minimal error in the solution. We first analyze the behavior of this system by deriving moment equations for a specifically-formulated system with concentration fluctuations. Using these tools, we ascertain that the error between a Dirac and Gaussian-based solution can be reduced to an arbitrary level, determined by choice of particle number and the parameters of the system. Further, we use these moment equation solutions to inform our criteria for matching the Dirac and Gaussian particle simulation solutions, essentially deducing the optimal kernel half-width from our choice of particle number by using one of three matching criteria.

The specific time matching criteria is the simplest method to implement. One need only choose the desired number of Gaussian particles, $N_G$, and the time at which solutions should be close, $t^*$, to compute the appropriate kernel size for the simulation. However, of the three methods, this one yielded the worst overall performance, as it is susceptible to error-inducing phenomena. First, the ratio of kernel half-width to domain size, $\ell_G/\Omega$, increases both as $N_G/N_\delta$ decreases and as $t^*$ increases. This inevitably leads to undesirable domain effects when $\ell_G\gtrapprox0.12$ in one dimension. Additionally, the particle solutions do not intersect as reliably near $t^*$ as the moment equation solutions do. The assumption that $\vert \bar C_\delta-\bar C_G\vert\approx0$ on the interval $[0,t^*]$, allowing for cancellation of the exponential integral terms in \eqref{g_delta_g}, seems not to hold up in the particle case, presumably due to the increasing effect of third-order moments as $N_G/N_\delta$ decreases.

The least squares matching criteria requires a numerical minimization of the error between moment equation solutions, but it provides a markedly better ability to minimize error in the particle simulations. This is due to the fact that, rather than assuming closeness of solutions for cancellation in order to minimize error at a single time, the algorithm minimizes the overall error in the Dirac and Gaussian moment equation solutions across a range of time values. As a result, it appears that the only source of error between our moment equation and particle tracking simulations is the neglect of third-order moments in the moment equations.

The emergent and persistent ``moat" behavior, discussed in Section \ref{sec:moats}, caused by the choice of a single kernel size for the duration of a simulation, led to an approach that allows the kernel half-width to grow with time. While we do not have the theoretical framework to fully analyze this method of matching, it shows strong potential as a method for matching Dirac and Gaussian solutions. One caveat, however, is that, while this variable $\ell_G$ method is able to minimize or eliminate domain effects at early time, it is still susceptible to domain effects at late time if $\ell_G(t)$ becomes larger than approximately 12\% of the domain.

To summarize, under the circumstances considered here, the least squares matching criteria is the most effective means for using up to an order of magnitude fewer Gaussian particles and attaining similar results to the corresponding Dirac particle solution. The single caveat is to keep the ratio of $\ell_G/\Omega$ smaller than approximately 0.12, in order to avoid error induced by domain effects.



\appendix
\section{Initial covariance relations for the particle tracking method}
\label{appa}

Following a similar approach to {\em Paster et al.} \cite{Paster_JCP}, we wish to derive representations for initial auto- and cross-covariance in concentration fluctuations for the particle tracking model. It is assumed that initial particle positions are independent and identically distributed and determined by drawing from a random uniform distribution. Throughout, we assume $\Omega$ is chosen to be sufficiently large, so as to exclude boundary effects caused by utilizing a finite computational domain (we investigate the effects of such finite domains numerically in Section \ref{sec:domain_effects}).

The concentrations composed of particles will be of the form
\beq
\bal
C_i(t,x)&=\sum_{k=1}^N\int_\Omega m_p\phi_i(x-z)\delta(z-x_k)dz,\ \ \ i=A,B,\\
&=m_p\sum_{k=1}^N\phi_i(x-x_k),
\eal
\eeq
where $x_k(t)$ is the position of the $k^{\text{th}}$ of $N$ particles, $m_p$ is the mass of a single particle (particle number and mass are assumed to be equal for $A$ and $B$ particles, for simplicity in calculations), and $\Omega$ is a $d$-dimensional domain. The kernel $\phi_i$ is defined to be symmetric ($\phi_i(x-z)=\phi_i(z-x)$), have units $L^{-d}$, and for $\Omega$ suitably large, to integrate approximately to unity $\int_\Omega\phi_i(x)dx\approx1$. For simplification in computation, we assume this approximation to be close enough so as to consider this integral equal to 1.

Thus, we may represent average concentration
\beq
\bal
\bar C_i(t)=E\left[C_i(t,x)\right]&=\frac{1}{\Omega}\int_\Omega C_i(t,x)dx\\
&=\frac{Nm_p}{\Omega}.
\eal
\eeq

Using the relation in \eqref{fluctuation}, initial autocovariance may be represented by $\bar{C_i'(0,x)C_i'(0,y)}=\bar{C_i(0,x)C_i(0,y)}-\bar C_i(0)^2$, yielding
\bln
\begin{alignat}{2}
\label{eq1}
\overline{C_i'(x)C_i'(y)}&=m_p^2\int\genfrac{}{}{0pt}{0}{\dots}{_{\Omega^{N}}}\int\Bigg[\sum_{k=1}^N\sum_{l=1}^N& &\phi_i(x-x_k)\phi_i(y-x_l)\\
& & &F(x_1,\dots,x_N)\Bigg]dx_1\dots dx_N-\left(\frac{Nm_p}{\Omega}\right)^2,\nonumber
\end{alignat}
\eln
where $F(x_1,\dots,x_N)$ is the joint pdf for the randomly distributed particles. Since the particle positions are independent of each other, $F(x_1,\dots,x_N)=F(x_1)\dots F(x_N)=\Omega^{-N}$, so the first term in (\ref{eq1}) becomes
\beq
\label{cross_tot_conc}
m_p^2\int\genfrac{}{}{0pt}{0}{\dots}{_{\Omega^{N}}}\int\left[\sum_{k=1}^N\sum_{l=1}^N\phi_i(x-x_k)\phi_i(y-x_l)\frac{1}{\Omega^{N}}\right]dx_1\dots dx_N.
\eeq

Now, first considering the cases above where $k\neq l$, we have

\beq
\bal
m_p^2\sum_{{\mathclap{\substack{k=1\\ k\neq l}}}}^N\sum_{{\mathclap{\substack{l=1\\ l\neq k}}}}^N\bigg[\int\genfrac{}{}{0pt}{0}{\dots}{_{\Omega^{N}}}\int&\left[\phi_i(x-x_k)\phi_i(y-x_l)\frac{1}{\Omega^N}\right]dx_1\dots dx_N\bigg]\\
&=\frac{N(N-1)m_p^2}{\Omega^{2}}\int_{\Omega}\phi_i(x-x_k)dx_k\int_{\Omega}\phi_i(y-x_l)dx_l\\
&=\frac{N(N-1)m_p^2}{\Omega^{2}},
\eal
\eeq
and when $k=l$
\beq
\frac{Nm_p^2}{\Omega}\int_\Omega\phi_i(x-x_k)\phi_i(y-x_k)dx_k=\frac{Nm_p^2}{\Omega}(\phi_i\star\phi_i)(x-y),
\eeq
where $\star$ denotes convolution. Thus, we have a general representation for initial autocovariance
\beq
\overline{C_i'(x)C_i'(y)}=\frac{Nm_p^2}{\Omega}\left[(\phi_i\star\phi_i)(x-y)-\frac{1}{\Omega}\right].
\eeq
For the specific Dirac delta and Gaussian kernel choices
\beq
\bal
\phi_\delta(x-y)&=\delta(x-y),\\
\phi_G(x-y)&=\frac{1}{(2\pi\ell_G^2)^{d/2}}e^{-\frac{\vert x-y\vert^2}{2\ell_G^2}},
\eal
\eeq
we have the following initial autocovariance structures
\beq
\label{spec_kernel_covs}
\begin{aligned}
\hat f_\delta(x-y)&=\frac{N_\delta m_\delta^2}{\Omega}\left[\delta(x-y)-\frac{1}{\Omega}\right]\\
&=C_0m_\delta\left[\delta(x-y)-\frac{1}{\Omega}\right],\\
\hat f_G(x-y)&=\frac{N_G m_G^2}{\Omega}\left[\frac{1}{(4\pi\ell_G^2)^{d/2}}e^{-\frac{\vert x-y\vert^2}{4\ell_G^2}}-\frac{1}{\Omega}\right]\\
&=C_0m_G\left[\frac{1}{(4\pi\ell_G^2)^{d/2}}e^{-\frac{\vert x-y\vert^2}{4\ell_G^2}}-\frac{1}{\Omega}\right].
\end{aligned}
\eeq

Next we calculate the initial cross-covariance, which is assumed to be reflexive as to particle type, so
\beq
\label{general_cross}
\bal
\bar{C'_A(0,x)C'_B(0,y)}&=\bar{C'_B(0,x)C'_A(0,y)}\\
&=\bar{C_A(0,x)C_B(0,y)}-\bar C_A(0)\bar C_B(0).
\eal
\eeq
Considering the first term in \eqref{general_cross}
\beq
\label{cross_calc_1}
\bar{C_A(x)C_B(y)}=\frac{m_p^2}{\Omega^{2N}}\int\genfrac{}{}{0pt}{0}{\dots}{_{\Omega^{2N}}}\int\sum_{k=1}^{N}\sum_{l=1}^{N}\phi_A(x-x_k)\phi_B(y-y_l)dx_1,\dots,dx_{N},dy_1,\dots,dy_{N}.
\eeq
In contrast to the autocovariance calculations, $x_k$ is never equal to $y_l$, since every $A$ and $B$ particle combination has probability zero of occupying the same position. As a result, we have
\beq
\frac{N^2m_p^2}{\Omega^{2}}\int_{\Omega}\int_{\Omega}\phi_A(x-x_k)\phi_B(y-y_l)dx_kdy_l=\frac{N^2m_p^2}{\Omega^{2}}.
\eeq
Now, considering the second term in \eqref{general_cross}, we find
\beq
\bar C_A\bar C_B=\frac{N^2m_p^2}{\Omega^{2}}.
\eeq
As such, it is apparent that our initial cross-covariance is identically equal to zero.



\section{Estimates for error propagation}
\label{appb}

In order justify the error estimates of Section \ref{sec:properror}, we first show some properties of the mean concentration equations.  Recall that, regardless of the choice of correlation structure, the mean concentration $\overline{C}$ satisfies
\beq
\label{CbarApp}
\left.
\begin{array}{l}
\displaystyle \frac{d\overline{C}}{dt} = -k \left [ \overline{C}^2 + g(t) \right ]\\
\\
\overline{C}(0) = C_0,
\end{array}
\right \}
\eeq
where
\beq
\label{geqn}
g(t) = \frac{1}{2(8\pi Dt)^{d/2}} \int f(0,z,x) e^{-\frac{\vert x - z \vert^2}{8Dt}} \ dz \left [-1 + \exp \left ( -4k\int_0^t \overline{C}(\tau) \ d\tau \right ) \right ].
\eeq

Since $g(t) = 0$ when $\overline{C} \equiv 0$, we see that $\overline{C} \equiv 0$ is an equilibrium point of the differential equation.  
Thus, because we take $C_0 > 0$, it follows that $\overline{C}(t) > 0$ for every $t \geq 0$ by uniqueness of solutions to (\ref{CbarApp}).
With this, we see that $g(t) < 0$ for every $t > 0$ because $\frac{1}{(8\pi Dt)^{d/2}} \int f(0,z,x) e^{-\frac{\vert x - z \vert^2}{8Dt}} \ dz \geq 0$ for sufficiently large $\Omega$ and the last term in \eqref{geqn} is negative.
Therefore, using the negativity of $g$, we find from (\ref{CbarApp})
\bln$$\frac{d\overline{C}}{dt} > -k\overline{C}^2,$$\eln
for all $t > 0$.
Because $\overline{C}(t) > 0$, we may use the separable nature of this inequality to find
\bln$$\frac{d}{dt} \left ( \overline{C}(t)^{-1} \right ) < k,$$\eln
which, upon integrating, finally implies for every $t > 0$
\bln$$ \overline{C}(t) > \frac{C_0}{1+kC_0 t}.$$\eln
Since both $\overline{C}_\delta$ and $\overline{C}_G$ satisfy a differential equation of this type, each must satisfy the same lower bound, namely,
 {\beq
\label{Clower}
\overline{C}_p(t) > \frac{C_0}{1+kC_0 t},
\eeq
for $p=\delta, G$.  This result merely displays the influence of the cross-covariance function so that the negative values of $g_p$ arising from spatial fluctuations cause the mean concentration to decrease at a slower rate than the well-mixed solution.

The above inequality further establishes an upper bound of growth on the exponential that appears in the $g_p$ terms.  In particular, using (\ref{Clower}) we find
\bln$$ -4k \int_0^t \overline{C}_p(\tau) d \tau < -4k \int_0^t \frac{C_0}{1+kC_0 \tau} d \tau = -4 \ln(1 + kC_0t),$$\eln
and thus
\beq
\label{expbound}
\exp \left ( -4k \int_0^t \overline{C}_p(\tau) d \tau \right ) < (1 + kC_0t)^{-4}.
\eeq


Finally, we consider the difference in forcing terms generated by different initial cross-covariances, namely
$\vert g_\delta(t) - g_G(t) \vert$
where each individual cross-covariance $g_p(t)$ is defined by (\ref{g}).
Subtracting the equations, we find
\bln
\begin{eqnarray*}
\vert g_\delta(t) - g_G(t) \vert & = & \left \vert \psi_\delta(t) \left (-1 + e^{-4k \int_0^t \overline{C}_\delta(\tau) d \tau} \right) - \psi_G(t) \left (-1 + e^{ -4k \int_0^t \overline{C}_G(\tau) d \tau} \right ) \right \vert\\
& \leq & I + II,
\end{eqnarray*}
\eln
where
\bln$$I :=  \left ( 1 - e^{-4k \int_0^t \overline{C}_G(\tau) d \tau} \right ) \vert \psi_\delta(t) - \psi_G(t) \vert ,$$\eln
and
\bln$$II := \vert \psi_\delta(t) \vert \left \vert e^{-4k \int_0^t \overline{C}_\delta(\tau) d \tau} - e^{-4k \int_0^t \overline{C}_G(\tau) d \tau} \right \vert.$$\eln
We note that near $t = 0$ each these terms tend to zero, and this can be verified by a simple use of L'Hospital's Rule.  Thus, there is no singularity at  $t = 0$ and we will concentrate on bounding these terms for $t$ suitably large.
Since the exponential in $I$ is nonnegative, we find
\bln
\begin{eqnarray*}
I & \leq & \vert \psi_\delta(t) -\psi_G(t) \vert\\
& \leq & \frac{1}{2}C_0 \left ( \frac{1}{\Omega} \vert m_\delta - m_G \vert + \left \vert \frac{m_\delta}{(8\pi Dt)^{d/2}} - \frac{m_G}{(4\pi \ell_G^2 + 8\pi Dt)^{d/2}} \right \vert \right )\\
& =: & I_A + I_B.
\end{eqnarray*}
\eln
Now, because $C_0\Omega = m_p N_p$ for $p=\delta, G$ the first term is exactly
$$I_A = \frac{1}{2}C_0^2 \left ( \frac{1}{N_G} - \frac{1}{N_\delta} \right ) = \frac{1}{2}C_0^2 \Delta.$$ 
Using the Mean Value Theorem on the function $h(x) = (x + 8\pi Dt)^{-d/2}$
so that
$$\vert h(0) - h(4\pi \ell_G^2)  \vert \leq 4\pi \ell_G^2 \max_{x \in [0,4\pi \ell_G^2]} \vert h'(x) \vert \leq \frac{2\pi d \ell_G^2}{(8\pi Dt)^{1+ d/2}} ,$$
the second term satisfies
\begin{eqnarray*}
I_B & \leq & \frac{1}{2} C_0 \left (\frac{\vert m_\delta - m_G \vert}{(8\pi Dt)^{d/2}} + m_G \left \vert \frac{1}{(8\pi Dt)^{d/2}} - \frac{1}{(4\pi \ell_G^2 + 8\pi Dt)^{d/2}} \right \vert \right )\\
& \leq & \frac{C_0^2 \Delta}{2\Omega (8\pi Dt)^{d/2}} + \frac{\pi dC_0 \ell_G^2m_G}{(8\pi Dt)^{1 + d/2}}.
\end{eqnarray*}

To estimate $II$, we use (\ref{psid}) to find
\bln$$\vert \psi_\delta(t) \vert \leq \frac{1}{2\Omega}C_0 m_\delta$$\eln
for $t \geq \frac{\Omega^{2/d}}{8\pi D}$.
Then, using (\ref{expbound}) it follows that
\bln$$II \leq \frac{C_0 m_\delta}{2\Omega} (1 + kC_0t)^{-4}.$$\eln
%
Recomposing these expressions within the inequality for the $g$ terms, we have
\bln$$\vert g_\delta(t) - g_G(t) \vert \leq  \frac{1}{2}C_0^2 \Delta + \frac{C_0^2 \Delta}{2\Omega (8\pi Dt)^{d/2}} + \frac{\pi dC_0 \ell_G^2m_G}{(8\pi Dt)^{1 + d/2}} +  \frac{C_0 m_\delta}{2\Omega}(1 + kC_0t)^{-4}.$$\eln
Finally, for simplicity we will assume within this analysis that $8\pi D > kC_0$ so that for $t \geq \frac{1}{8\pi D - kC_0}$, we have the inequality
$$\frac{1}{8 \pi Dt} \leq \frac{1}{1 + kC_0t},$$
and thus
\beq
\label{gbound}
\vert g_\delta(t) - g_G(t) \vert \leq  \frac{1}{2}C_0^2 \Delta + \frac{C_0^2 \Delta}{2\Omega (1 + kC_0t)^{d/2}} + \frac{\pi dC_0 \ell_G^2m_G}{(1 + kC_0t)^{1+ d/2}} +  \frac{C_0 m_\delta}{2\Omega(1 + kC_0t)^{4}}.
\eeq

It should be noted that a nearly identical analysis can be conducted when $8\pi D \leq kC_0$, which merely results in replacing the $(1 + kC_0t)$ terms in \eqref{gbound} with $8 \pi Dt$, which is also $\mathcal{O}(t)$ as $t \to \infty$.  As can be seen in Section $4$, this will yield identical conclusions and estimates on the error incurred by using the kernel-based method rather than the $\delta$-based reactive particle tracking method.


\section*{References}

\bibliography{reaction_proposal}

\begin{thebibliography}{10}
\expandafter\ifx\csname url\endcsname\relax
  \def\url#1{\texttt{#1}}\fi
\expandafter\ifx\csname urlprefix\endcsname\relax\def\urlprefix{URL }\fi
\expandafter\ifx\csname href\endcsname\relax
  \def\href#1#2{#2} \def\path#1{#1}\fi

\bibitem{battiato1}
I.~Battiato, D.~M. Tartakovsky, Applicability regimes for macroscopic models of
  reactive transport in porous media, J. Contam. Hydrol. 120-121 (2011) 18--26.
\newblock \href {http://dx.doi.org/10.1016/j.jconhyd.2010.05.005}
  {\path{doi:10.1016/j.jconhyd.2010.05.005}}.

\bibitem{battiato2}
I.~Battiato, D.~M. Tartakovsky, A.~M. Tartakovsky, T.~Scheibe, On breakdown of
  macroscopic models of mixing-controlled heterogeneous reactions in porous
  media, Adv. Water Resour. 32 (2009) 1664--1673.
\newblock \href {http://dx.doi.org/10.1016/j.advwatres.2009.08.008}
  {\path{doi:10.1016/j.advwatres.2009.08.008}}.

\bibitem{Schwede_sample_pdf}
R.~L. Schwede, O.~A. Cirpka, W.~Nowak, I.~Neuweiler,
  \href{http://dx.doi.org/10.1029/2007WR006668}{Impact of sampling volume on
  the probability density function of steady state concentration}, Water
  Resources Research 44~(12) (2008) W12433, w12433.
\newblock \href {http://dx.doi.org/10.1029/2007WR006668}
  {\path{doi:10.1029/2007WR006668}}.
\newline\urlprefix\url{http://dx.doi.org/10.1029/2007WR006668}

\bibitem{Chiogna2013}
G.~Chiogna, A.~Bellin, \href{http://dx.doi.org/10.1002/wrcr.20200}{Analytical
  solution for reactive solute transport considering incomplete mixing within a
  reference elementary volume}, Water Resources Research 49~(5) (2013)
  2589--2600.
\newblock \href {http://dx.doi.org/10.1002/wrcr.20200}
  {\path{doi:10.1002/wrcr.20200}}.
\newline\urlprefix\url{http://dx.doi.org/10.1002/wrcr.20200}

\bibitem{Sweby1984}
P.~K. Sweby, \href{http://www.jstor.org/stable/2156939}{High resolution schemes
  using flux limiters for hyperbolic conservation laws}, SIAM Journal on
  Numerical Analysis 21~(5) (1984) 995--1011.
\newline\urlprefix\url{http://www.jstor.org/stable/2156939}

\bibitem{Leonard1991}
B.~Leonard,
  \href{http://www.sciencedirect.com/science/article/pii/004578259190232U}{The
  ultimate conservative difference scheme applied to unsteady one-dimensional
  advection}, Computer Methods in Applied Mechanics and Engineering 88~(1)
  (1991) 17--74.
\newline\urlprefix\url{http://www.sciencedirect.com/science/article/pii/004578259190232U}

\bibitem{Leveque2002}
R.~J. Le{V}eque, Finite Volume Methods for Hyperbolic Problems, Cambridge
  University Press, 2002.

\bibitem{Bokanowski2007}
O.~Bokanowski, H.~Zidani,
  \href{http://dx.doi.org/10.1007/s10915-005-9017-0}{Anti-dissipative schemes
  for advection and application to {H}amilton--{J}acobi--{B}ellmann equations}
  30~(1) (2007) 1--33--.
\newline\urlprefix\url{http://dx.doi.org/10.1007/s10915-005-9017-0}

\bibitem{Gillespie1977}
D.~T. Gillespie, Exact stochastic simulation of coupled chemical reactions, J.
  Phys. Chem. 81~(25) (1977) 2340--2361.

\bibitem{Gillespie2013}
D.~T. Gillespie, A.~Hellander, L.~R. Petzold,
  \href{http://scitation.aip.org/content/aip/journal/jcp/138/17/10.1063/1.4801941}{Perspective:
  Stochastic algorithms for chemical kinetics}, The Journal of Chemical Physics
  138~(17) (2013) 170901.
\newblock \href {http://dx.doi.org/http://dx.doi.org/10.1063/1.4801941}
  {\path{doi:http://dx.doi.org/10.1063/1.4801941}}.
\newline\urlprefix\url{http://scitation.aip.org/content/aip/journal/jcp/138/17/10.1063/1.4801941}

\bibitem{Isaacson}
S.~A. Isaacson,
  \href{http://scitation.aip.org/content/aip/journal/jcp/139/5/10.1063/1.4816377}{A
  convergent reaction-diffusion master equation}, The Journal of Chemical
  Physics 139~(5) (2013) 054101.
\newblock \href {http://dx.doi.org/http://dx.doi.org/10.1063/1.4816377}
  {\path{doi:http://dx.doi.org/10.1063/1.4816377}}.
\newline\urlprefix\url{http://scitation.aip.org/content/aip/journal/jcp/139/5/10.1063/1.4816377}

\bibitem{HellanderPRE}
S.~Hellander, L.~Petzold,
  \href{http://link.aps.org/doi/10.1103/PhysRevE.93.013307}{Reaction rates for
  a generalized reaction-diffusion master equation}, Phys. Rev. E 93 (2016)
  013307.
\newblock \href {http://dx.doi.org/10.1103/PhysRevE.93.013307}
  {\path{doi:10.1103/PhysRevE.93.013307}}.
\newline\urlprefix\url{http://link.aps.org/doi/10.1103/PhysRevE.93.013307}

\bibitem{Benson_react}
D.~A. Benson, M.~M. Meerschaert,
  \href{http://dx.doi.org/10.1029/2008WR007111}{Simulation of chemical reaction
  via particle tracking: Diffusion-limited versus thermodynamic rate-limited
  regimes}, Water Resour. Res. 44 (2008) W12201.
\newblock \href {http://dx.doi.org/10.1029/2008WR007111}
  {\path{doi:10.1029/2008WR007111}}.
\newline\urlprefix\url{http://dx.doi.org/10.1029/2008WR007111}

\bibitem{vanZon}
J.~S. van Zon, P.~R. ten Wolde,
  \href{http://link.aps.org/doi/10.1103/PhysRevLett.94.128103}{Simulating
  biochemical networks at the particle level and in time and space: Green's
  function reaction dynamics}, Phys. Rev. Lett. 94 (2005) 128103.
\newblock \href {http://dx.doi.org/10.1103/PhysRevLett.94.128103}
  {\path{doi:10.1103/PhysRevLett.94.128103}}.
\newline\urlprefix\url{http://link.aps.org/doi/10.1103/PhysRevLett.94.128103}

\bibitem{Bolster2012}
D.~Bolster, P.~de~Anna, D.~A. Benson, A.~M. Tartakovsky,
  \href{http://www.sciencedirect.com/science/article/pii/S030917081100217X}{Incomplete
  mixing and reactions with fractional dispersion}, Advances in Water Resources
  37~(0) (2012) 86--93.
\newblock \href {http://dx.doi.org/10.1016/j.advwatres.2011.11.005}
  {\path{doi:10.1016/j.advwatres.2011.11.005}}.
\newline\urlprefix\url{http://www.sciencedirect.com/science/article/pii/S030917081100217X}

\bibitem{Bolster_mass}
D.~Bolster, A.~Paster, D.~A. Benson,
  \href{http://dx.doi.org/10.1002/2015WR018310}{A particle number conserving
  {L}agrangian method for mixing-driven reactive transport}, Water Resources
  Research 52~(2) (2016) 1518--1527.
\newblock \href {http://dx.doi.org/10.1002/2015WR018310}
  {\path{doi:10.1002/2015WR018310}}.
\newline\urlprefix\url{http://dx.doi.org/10.1002/2015WR018310}

\bibitem{Benson_arbitrary}
D.~A. Benson, D.~Bolster, Arbitrarily complex chemical reactions on particles,
  Water Resources Research Submitted.

\bibitem{Paster_JCP}
A.~Paster, D.~Bolster, D.~A. Benson,
  \href{http://www.sciencedirect.com/science/article/pii/S0021999114000473}{Connecting
  the dots: {S}emi-analytical and random walk numerical solutions of the
  diffusion--reaction equation with stochastic initial conditions}, Journal of
  Computational Physics 263~(0) (2014) 91--112.
\newblock \href {http://dx.doi.org/10.1016/j.jcp.2014.01.020}
  {\path{doi:10.1016/j.jcp.2014.01.020}}.
\newline\urlprefix\url{http://www.sciencedirect.com/science/article/pii/S0021999114000473}

\bibitem{Dani_kernel}
M.~Rahbaralam, D.~Fern{\`a}ndez-Garcia, X.~Sanchez-Vila,
  \href{http://www.sciencedirect.com/science/article/pii/S0021999115006233}{Do
  we really need a large number of particles to simulate bimolecular reactive
  transport with random walk methods? {A} kernel density estimation approach},
  Journal of Computational Physics 303 (2015) 95 -- 104.
\newblock \href {http://dx.doi.org/http://dx.doi.org/10.1016/j.jcp.2015.09.030}
  {\path{doi:http://dx.doi.org/10.1016/j.jcp.2015.09.030}}.
\newline\urlprefix\url{http://www.sciencedirect.com/science/article/pii/S0021999115006233}

\bibitem{kang_prl}
K.~Kang, S.~Redner, Scaling approach for the kinetics of recombination
  processes, Physical Review Letters 52 (1984) 955--958.

\bibitem{Ovchinnikov1978}
A.~A. Ovchinnikov, Y.~B. Zeldovich, Role of density fluctuations in bimolecular
  reaction kinetics, Chem. Phys. 28 (1978) 215--218.

\bibitem{KANG1984}
K.~Kang, S.~Redner, Fluctuation effects in {S}moluchowski reaction-kinetics,
  Physical Review A 30~(5) (1984) 2833--2836.
\newblock \href {http://dx.doi.org/10.1103/PhysRevA.30.2833}
  {\path{doi:10.1103/PhysRevA.30.2833}}.

\bibitem{Toussaint}
D.~Toussaint, F.~Wilczek, Particle--antiparticle annihilation in diffusive
  motion, J. Chem. Phys. 78~(5) (1983) 2642--2647.

\bibitem{wrr_incomplete}
A.~Tartakovsky, P.~de~Anna, T.~L. Borgne, A.~Balter, D.~Bolster, Effect of
  spatial concentration fluctuations on effective kinetics in
  diffusion-reaction systems,, Water Resources Research 48 (2012) W02526.

\bibitem{Paster_WRR}
A.~Paster, D.~Bolster, D.~A. Benson, Particle tracking and the
  diffusion-reaction equation, Water Resour. Res. 49 (2013) 1--6.
\newblock \href {http://dx.doi.org/10.1029/2012WR012444}
  {\path{doi:10.1029/2012WR012444}}.

\end{thebibliography}

\end{document}